\documentclass[11pt, oneside]{amsart}

\usepackage{multirow}
\usepackage{amsmath,amscd,amssymb}
\usepackage{xspace}
\usepackage{cite,alltt}
\usepackage{enumitem}
\setlist[itemize]{topsep=2pt,parsep=0pt,partopsep=2pt}
\setlist[enumerate]{topsep=2pt,parsep=0pt,partopsep=2pt}

\usepackage{comment}
\usepackage{mathtools, color}
\usepackage{multirow}
\usepackage{rotating}
\usepackage{subfigure}
\usepackage{epsfig}
\usepackage{url}
\usepackage{graphicx}
\usepackage{stackengine}
\usepackage{mathdots}
\usepackage[english]{babel}
\usepackage{csquotes}
\usepackage{mathrsfs}
\usepackage{overpic}
\usepackage{algorithmic} 
\usepackage{algorithm}
\usepackage{blindtext}
\usepackage{mathrsfs}
\usepackage{hyperref}
\usepackage{tikz}
\usepackage{pgfplots}
\usetikzlibrary{positioning}
\usepackage[export]{adjustbox}

\usepackage{bbm}
\usepackage{bm}
\usepackage{mathrsfs}
\usepackage{arydshln}
\usepackage{url,cite}
\usepackage{enumerate}
\usepackage[title]{appendix}
\usepackage{color}

\usepackage[normalem]{ulem}
\usepackage{soul}

\definecolor{g1}{HTML}{2da02d}
\definecolor{r1}{HTML}{d62728}
\definecolor{b1}{HTML}{2178b5}

\definecolor{Cerulean}{HTML}{00A2E3}
\definecolor{JungleGreen}{HTML}{00A99A}
\definecolor{Orange}{HTML}{F58137}

\newtheorem{remark}{Remark}
\usepackage{multirow,booktabs}
\usepackage{makecell}
\usepackage{subfigure}
\numberwithin{equation}{section}

\newenvironment{fequation}{\begin{equation}}{\end{equation}}
\newtheorem{theorem}{Theorem}[section]

\newtheorem{proposition}[theorem]{Proposition}
\newtheorem{corollary}[theorem]{Corollary}
\newtheorem{lemma}[theorem]{Lemma}
\newtheorem{assumption}{Assumption}
\newtheorem{definition}[theorem]{Definition}
\def\R{\mathbb{R}}

\def\S{\mathcal{S}}
\def\Sph{\mathbb{S}}

\def\d{\mathrm{d}}

\def\qed{~\relax\ifmmode\hskip2em \Box
	\else\unskip\nobreak\hskip1em \hfill$\Box$
	\fi \newline}

\newmuskip\pFqmuskip

\newcommand*\pFq[6][8]{%
	\begingroup 
	\pFqmuskip=#1mu\relax
	\mathchardef\normalcomma=\mathcode`,
	\mathcode`\,=\string"8000
	\begingroup\lccode`\~=`\,
	\lowercase{\endgroup\let~}\pFqcomma
	{}_{#2}F_{#3}{\left(\genfrac..{0pt}{}{#4}{#5};#6\right)}%
	\endgroup
}
\newcommand{\pFqcomma}{{\normalcomma}\mskip\pFqmuskip}

\newcommand{\hgeom}[2]{\,{}_{#1\!}F_{#2}}
\def\R{\mathbb{R}}
\def\S{\mathbb{S}}

\def\P{\mathcal{P}}

\def\d{\mathrm{d}}
\def\i{\mathrm{i}}

\def\IntSph{\int_{\S^d}}
\def\ConvPhi{\mathcal{C}_{\varphi_{\rho}}}
\def\Nro{N_{\rho}}

\def\QuadSet{\mathcal{Q}_{X}}
\def\scaleKer{\varphi_{\rho}}
\newcommand{\NN}{\mathbb{N}}

\title[Spherical quasi-interpolation]{Scaled zonal kernel quasi-interpolation on spheres}
\author{Zhengjie Sun}
\address{School of Mathematics and Statistics,  Nanjing University of Science and Technology, Nanjing, China}
\email{zhengjiesun@njust.edu.cn}

\author{Wenwu Gao}
\address{Corresponding author. School of Big Data and Statistics, Anhui University, Hefei, China}
\email{wenwugao528@163.com}

\author{Xingping Sun}
\address{Department of Mathematics, Missouri State University, Springfield, MO 65897, USA}
\email{xsun@missouristate.edu}

\graphicspath{{Figures_PeriKer/}}

\usepackage{geometry}
\begin{document}
	\include{Data_sphereQI}
	\begin{abstract}
		We propose and study a new quasi-interpolation method on spheres featuring the following two-phase construction and analysis. In Phase I, we  analyze and characterize a large family of zonal kernels   
		(e.g., the spherical version of Poisson kernel, Gaussian, compactly-supported radial kernels),
		so that the underlying spherical convolution operators (upon the introduction of a scaling parameter) attain a high-order approximation to target functions. 
		In Phase II, we discretize the spherical integrals utilizing quadrature rules to produce the final quasi-interpolants. Numerical experiments demonstrate that the new approximation algorithm 
		is robust and amenable to integrated as well as distributed ways of implementation. Moreover, the underlying error-analysis shows that by fine-tuning the scaling parameter in the radial kernels employed, the resulting quasi-interpolants achieve a well-balanced trade-off between approximation and sampling errors. 
	\end{abstract}
	
	\maketitle
	
	\section{Introduction}
	Function approximation on spherical domains is a fundamental and important topic in various scientific fields, with applications ranging from cluster computing to data analysis and machine learning. 
	From a historical perspective, 
	polynomials  have been the earliest and perhaps the most-studied tool for function approximation \cite{buhmann2003radial-book}, \cite{wendland2004scattered}. 
	Despite its rich mathematical heritage and theoretical superiority, applicable limitations of polynomial approximation persist, which are reflected on the one hand by the inherent nature of polynomials (frequent  oscillation of high-degree polynomials) and on the other hand by construction methods of polynomial approximation; see \cite{sloan_2012Geom_filtered}.
	
	To facilitate applications of polynomial approximation, Sloan \cite{sloan_1995JAT_polynomial} introduced a  novel 
	quasi-interpolation method dubbed ``hyperinterpolation", the spherical version of which involves discretizing the following spherical convolution integral
	\[
	\int_{\S^d}f(y)D_n(x, y)\d\mu(y),
	\]
	where $D_n$ is the Dirichlet kernel associated with the $n$th partial sum of the Fourier-Legendre series of a target function $f \in C(\S^d),$ and $\mu$ the restriction to $\S^d$ of the Lebesgue measure on $\R^{d+1}.$  Sloan later \cite{sloan2011polynomial} 
	proposed ``filtered hyperinterpolation" in which kernels with locality properties are employed 
	to temper the ``hard" truncation of Fourier-Legendre series of a target function induced by the Dirichlet kernel. Sloan's hyperinterpolation method has attracted the attention of fellow researchers. Generalizations and improvements of it soon ensued. With the intention to handle massive spherical noisy data, Lin, Wang and Zhou \cite{lin_SINUM2021_distributed} utilized filtered hyperinterpolation in conjunction with a distributed learning technique to great effect. More recently, Mont\'{u}far and Wang  \cite{montufar_2022FoCom_distributed} provided a comprehensive study  of distributed filtered hyperinterpolation on manifolds.
	
	Wang and Sloan \cite{wang2017filtered} provided a   \textbf{two-phase} viewpoint of filtered hyperinterpolation \cite{sloan_2012Geom_filtered}: first filtered polynomial approximation \cite{sloan2011polynomial} and then  discretization of remaining Fourier-Legendre coefficients  via
	quadrature rules. This has allowed them to derive  error estimates of  filtered hyperinterpolation respectively in terms of the filtered polynomial approximation error and the discretization error.  Here we take the liberty of using the notations (to be introduced) in the beginning part of Section 2 and paraphrase 
	the spherical version of the filtered hyperinterpolation method as follows. To begin with, express an $ f \in L_1(\S^d)$ by its Fourier-Legendre series:
	\begin{equation}\label{Fourier-coefficients}
		f=\sum_{\ell=0}^{\infty}\sum_{k=1}^{N(d,\ell)}\langle f,Y_{\ell, k}\rangle Y_{\ell, k}, \quad {\rm where} \quad \langle f,Y_{\ell, k}\rangle:= \int_{\S^d} f(x) Y_{\ell, k}(x)\d \mu(x).     
	\end{equation}
	This is followed by applying the 
	multiplier operator m$=\left\{G_a\Big(\frac{\ell}{L}\Big)\right\}_{\ell=0}^\infty$ ($L \in \NN$) to the series above:
	\begin{equation}\label{multiplier-action}
		\mbox{m}(f)=\sum_{\ell=0}^{\infty}G_a\Big(\frac{\ell}{L}\Big)\sum_{k=1}^{N(d,\ell)}\langle f,Y_{\ell, k}\rangle Y_{\ell, k}.    
	\end{equation}
	Here the function  $G_a \in C^\kappa[0,\infty)$ is referred to as a ``filter", which is compactly-supported on the interval $[0, 2a]$,  ($1<a<\infty$ being an adjustable parameter) and satisfies $G_a(x)=1$, for $x\in [0,a]$, and $G_a(x)=0$, for $x\geq 2a$. The degree of smoothness $\kappa$ for the filter is  chosen appropriately so that the resulting multiplier operator is bounded on some Sobolev spaces. 
	For more details on specifics of the filter, we refer readers to  \cite{sloan2011polynomial}, \cite{sloan_2012Geom_filtered}, \cite{wang2017filtered} and the references therein. Equivalently, Equation \eqref{multiplier-action} can be written in the following spherical convolution form ((2.2) in \cite{wang2017filtered}):
	\begin{equation}\label{spherical-convolution}
		\mathcal{V}_{L,G_a}(f)(x):=\int_{\S^d}f(y)\Phi_{L,G_a}(x, y)\d\mu(y),   
	\end{equation}  in which
	\begin{equation}\label{hyperkernel}
		\Phi_{L,G_a}(x, y) :=\sum_{\ell=0}^{\infty}G_a\Big(\frac{\ell}{L}\Big)K_{\ell}(x,y)
	\end{equation}
	with  $K_{\ell}(x,y)=\sum_{k=1}^{N(d,\ell)} Y_{\ell, k}(x)Y_{\ell, k}(y)$.  
	In the final stage, the integral on the right hand side of \eqref{spherical-convolution} is discretized using an $N$-point quadrature rule (of suitably high polynomial degree of precision) to give rise to the quasi-interpolant $\mathcal{V}_{L,G_a}^N(f)$ for the underlying filtered hyperinterpolation scheme:
	$$\mathcal{V}_{L,G_a}^N(f)(x)=\sum_{j=1}^Nf(x_j)\Phi_{L,G_a}(x, x_j)\omega_j$$ for some positive weights $\omega_j$, $j=1,2,\cdots, N$. 
	
	Filtered hyperinterpolation has been extensively studied in the literature; see  \cite{an_SINUM2012_regularized,hansen_IMA2009_OntheNorm,lin_SINUM2021_distributed,sloan2011polynomial,sloan_2012Geom_filtered,wang_2017ACHA_fully}  and references therein.  However, polynomial zonal-kernels have encountered computational challenges. In particular, the oscillatory nature of high-degree polynomials causes the underlying algorithm to be susceptible to noisy data.  More unsettling for some field scientists is the lack of clear ways in which  the positivity of filtered hyperinterpolation kernels shown in \eqref{hyperkernel} is determined. Except for the well-known Askey-Gasper inequality \cite{Askey-1995AmerMath-positive} which asserts that
	\[
	\sum_{j=0}^{n}\frac{C_{j}^{\alpha }(x)}{{2\alpha +j-1 \choose j}}\geq 0\qquad (x\geq -1,\,\alpha \geq 1/4),
	\]
	where $C_{j}^{\alpha }(x)$ denotes the Gegenbauer polynomials, not much is known about the positivity of the likes of filtered hyperinterpolation kernels. This concern has previously been voiced by Sloan and Womersley \cite{sloan_2012Geom_filtered}.
	
	In this paper, we propose and study an alternative to filtered hyperinterpolation. Let $\phi$ be a radial function on $\R^{d+1}$, that means that there is a function $\zeta \in C[0, \infty)$ such that $\phi(x)=\zeta(\|x\|)$ , where $\|\cdot\|$ is Euclidean norm.
	For $0 < \rho < 1,$  define
	$\phi_\rho(x)=
	\phi(\rho^{-1}x).$ Restrict the kernel $\phi_\rho(x-y)$
	on $\R^{d+1} \times \R^{d+1} $ to $\S^d \times \S^d$ subject to the law of cosine:
	\[
	\|x-y\|^2 = 2 - 2 x \cdot y, \quad x,y \in \S^d.
	\]
	By Funk-Hecke formula, we have
	\[
	\Lambda_{\phi,\rho}=\int_{\S^d} \phi_\rho(x-y) \d\mu(y) = c(d,\rho) \int^1_{-1} \phi_\rho(t) (1-t^2)^{(d-1)/2} \d t,
	\]
	where $c(d,\rho)$ is a constant depending only on $\rho$ and $d.$ Suppose $\Lambda_{\phi,\rho} \ne 0$ and that $\phi$ decays appropriately away from zero.
	Then we form an ``approximate identity" kernel $\varphi_\rho$ 
    by setting
	\begin{equation}\label{eq:DefiSphKernel}
		\varphi_\rho(x\cdot y)=(\Lambda_{\phi, \rho})^{-1} \phi_\rho(x -y),
	\end{equation}
	with which we construct a spherical convolution operator
	$$\ConvPhi(f)(x):=\int_{\S^d}f(y)\varphi_{\rho}(x\cdot y)\d\mu(y), \quad f \in L_p(\S^d).$$
	This general construction framework admits a large family of zonal kernels for approximation purposes. In the current paper, we investigate 
	several radial kernels restricted to the sphere, including Poisson kernel, Gaussian kernel, and compactly-supported radial kernels \cite{kuo2017multivariate}. Emphatically, we mention the spherical version of the classical Poisson kernel, which we will further elaborate in Section 4. We provide easily verifiable 
	conditions (see Assumption \ref{Assump1:Kernel}) under which we employ spherical convolution operators $\ConvPhi(f)$ to approximate target functions and gauge approximation orders in terms of  $\rho$ (see Theorem \ref{thm:SphConvConvergence}).  
	

 In implementing the approximation method, we 
 discretize the convolution integral using function values obtained on some finite subsets of $\S^d$. There exist in the numerical analysis literature numerous ways in which the convolution integrals are discretized to yield superior approximation powers. These include spherical 
$t$-designs, quasi-Monte Carlo methods, crosslet sparse grids, lattice rules, and stratified random samplings \cite{bourgain-1, bourgain-2, bourgain-3, buhmann-2021Adv-discretization,dick--2013ActaNumer--high, gao-2022NA-multivariate,kuo2017multilevel,kuo2021function,mhaskar-2001Mcom-spherical}. Here we choose the spherical quadrature rules studied in \cite{dai_2006JFA_multivariate,dai2013approximation,brauchart2007ConsAppr_numerical,hesse2010numerical-book} as a pilot example. This choice aligns with the methodology commonly employed in spherical hyperinterpolation \cite{sloan_1995JAT_polynomial,womersley_sloan_spherepoints,sloan_2012Geom_filtered} and offers a direct comparison between our new method and the conventional hyperinterpolation method. Furthermore, our computer simulations show that
in the context of noisy data our new method outperforms the conventional hyperinterpolation method;
see Figures \ref{fig.Compar_noisydata} and \ref{fig.Compar_CpuTime} in Section \ref{sec:NumerExamp}.

Upon discretizing  $C_{\varphi_{\rho}}(f)$ using an appropriately-chosen $|X|$-point spherical quadrature rule $(x_j, w_j)_{j=1}^{|X|}$, we
derive the quasi-interpolant: \begin{equation*}\label{eq:QIgenearlform}   
		\ConvPhi^X(f)(x):=\sum_{j=1}^{|X|}f(x_j)\varphi_{\rho}(x\cdot x_j)w_j.
	\end{equation*}

The literature of quasi-interpolation boasts a wide range of methods and tools and has been extensively  studied and broadly applied in diverse scientific fields, including surface reconstruction \cite{liu-wang2012quasi}, image processing \cite{costarelli2020comparison}, data fitting \cite{buhmann2003radial-book}, solving partial differential equations \cite{sun-jcp2019multi}. 
Special optimality and regularization properties pertaining to quasi-interpolation have also been discussed in the literature; see
\cite{gao2020optimality}. 
For quasi-interpolation on spheres, several existing approaches have been studied and collected in \cite{buhmann2022quasi-book}. These include tensor product trigonometric splines on spherical coordinates \cite{ganesh2006quadrature}, truncating Fourier coefficients \cite{gomes2001approximation}, and spline-based methods through sphere triangulation \cite{ibanez2010construction}.
However, the mesh-free and easy-implementation features of our method make it more applicable in a host of real world problems.
Furthermore, our method has the following two added advantages. (i) It does not entail computing Legendre polynomials of various degrees and orders when the  dimension of the ambient space changes. Our numerical experiments show that this enhances the stability of the underlying algorithms in various computing environments; (ii) Catering to practical needs, one may choose Taylor-made radial kernels to achieve desirable effects of quasi-interpolants, such as  positivity, monotonicity, and divergence/curl-free properties \cite{gao-fasshauer-fisher2022divergence}.

The paper is organized as follows.
In Section \ref{sec:Prelim}, we present preliminary materials, including the notations and definitions used throughout the paper. Section \ref{sec:sphQI} is dedicated to the construction of spherical convolution operators using zonal kernels. We provide a comprehensive characterization of zonal kernels that ensure the spherical convolution sequence converges to target functions with high-order accuracy. Additionally, we present a collection of important examples that serve as convolution approximation kernels and describe a construction technique for deriving high-order kernels from commonly used kernels. In Section \ref{sec:sphqi_scheme}, we discretize the spherical integrals using a specific quadrature rule and conclude the final quasi-interpolation scheme along with detailed error estimates. Finally, in Section \ref{sec:NumerExamp}, we present numerical results to validate the effectiveness of the proposed spherical quasi-interpolation method.


\section{Preliminaries}
\label{sec:Prelim}

We work on the unit sphere $\S^d:=\{x\in\R^{d+1}:\|x\|_2=1\}\subseteq \R^{d+1}$. The distance on the sphere is measured by the geodesic distance given by $\text{dist}(x,y)=\arccos(x\cdot y)$.  Let $\Delta_{\S}$ be the Laplace-Beltrami operator on the sphere. It is well known that the eigenvalues corresponding to the operator $-\Delta_{\S}$ are
$\lambda_{\ell}=\ell(\ell+d-1),~~\ell\in\mathbb{N}_0$.
The corresponding eigenfunctions are (spherical) polynomials called spherical harmonics, $\{Y_{\ell,k}\}_{k=1}^{N(d,\ell)}$ of degree $\ell$ and order $k$. Here, $N(d,\ell)$ is the dimension of this space, given by
\begin{equation}\label{eq:DimSphHarm}
	N(d,0)=1;~~N(d,\ell)=\frac{(2\ell+d-1)\Gamma(\ell+d-1)}{\Gamma(\ell+1)\Gamma(d)},~\ell\geq 1,
\end{equation}
with $N(d,\ell)=\mathcal{O}(\ell^{d-1})$ as $\ell\rightarrow\infty$. Moreover, the collection of all the $Y_{\ell,k}$'s forms an orthonormal basis for $L_2(\mathbb{S}^d)$. 
The space of all spherical polynomials $\Pi_n(\mathbb{S}^d)$ of degree $n$ or less is denoted by
\begin{equation*}
	\Pi_n(\mathbb{S}^d)=\text{span}\{Y_{\ell,k}:0\leq\ell\leq n, ~1\leq k\leq N(d,\ell)\}.
\end{equation*}
We employ the standard inner product for functions $f$, $g$ in  $L_2(\S^d)$,  
\begin{equation*}
	\langle f, g\rangle=\int_{\mathbb{S}^d}f(x)g(x)\d\mu(x),
\end{equation*}
where $\d\mu$ is the volume element of $\S^d$. For any function $f\in L_2(\mathbb{S}^d)$, its associated Fourier series converge in $L_2(\mathbb{S}^d)$,
\begin{equation}\label{eq:fexpansion}
	f=\sum_{\ell=0}^{\infty}\sum_{k=1}^{N(d,\ell)}\widehat{f}_{\ell,k}Y_{\ell,k},~~\widehat{f}_{\ell,k}=\int_{\S^d} f Y_{\ell,k}\d\mu.
\end{equation}
For $\sigma \ge 0$, the Sobolev space of order $\sigma$ on $\S^d$,  $H^{\sigma}(\S^d)$ is defined by 
$$H^{\sigma}(\S^d):=\{f\in L_2(\S^d):\|f\|_{H^{\sigma}(\S^d)}<\infty\},$$
where the norm is given by
\begin{equation}\label{eq:NormDef}
	\|f\|_{H^{\sigma}(\mathbb{S}^d)}^2=\sum_{\ell=0}^{\infty}\sum_{k=1}^{N(d,\ell)}(1+\ell)^{2\sigma}|\widehat{f}_{\ell,k}|^2,
\end{equation}
which is induced by the inner product:
$$\langle f,g\rangle = \sum_{\ell=0}^{\infty}\sum_{k=1}^{N(d,\ell)}(1+\ell)^{2\sigma}\widehat{f}_{\ell,k}\widehat{g}_{\ell,k}.$$
In particular, $H^{0}(\S^d)=L_2(\S^d).$
To construct spherical convolution operators, we will use zonal kernels on $\S^d\times \S^d,$ which admit the following symmetric representation in terms of spherical harmonics:
\begin{equation}\label{eq:zonalfunctions}
	\Psi(x,y)=\varphi(x\cdot y)=\sum_{\ell=0}^{\infty} a_{\ell}P_{\ell}(d+1;x\cdot y), 
\end{equation}
where $\varphi:[-1,1]\rightarrow\mathbb{R}$ is a continuous function and $P_{\ell}(d+1;t)$  the $(d+1)$-dimensional Legendre polynomial of degree $\ell$, normalized such that $P_{\ell}(d+1;1)=1$. 
The series on the right hand side of \eqref{eq:zonalfunctions} is referred to as the ``Fourier-Legendre expansion" of the zonal kernel $\varphi$, and $a_\ell$ the ``Fourier-Legendre coefficients". 
Convergence of the series is in the sense of Schwartz class distributions. 
An efficient way of calculating the Fourier-Legendre coefficients is via Funk-Hecke formula, which states that for every spherical harmonic of degree $\ell$ and order $k$, it holds that
\begin{equation}\label{eq:FunkHeckeFormula}
	\int_{\mathbb{S}^d}\varphi(x\cdot y)Y_{\ell,k}(y)\d\mu(y)=\widehat{\varphi}(\ell)Y_{\ell,k}(x),~
	\text{with} ~~\widehat{\varphi}(\ell)=\int_{-1}^1\varphi(t)P_{\ell}(d+1;t)\d t.
\end{equation}
This together with the addition formula for spherical harmonics (e.g. \cite[Page 10]{muller1966SphHarm})
leads to
\begin{equation}\label{kernel-coefficient}
	\varphi(x\cdot y)=\sum_{\ell=0}^{\infty}\widehat{\varphi}(\ell)K_{\ell}(x,y),
\end{equation}
in which $\widehat{\varphi}(\ell)=\frac{\omega_d}{N(d,\ell)}a_{\ell}$, where $\omega_{d}$ is the surface area of $\S^d$, and \begin{equation*}\label{eq:AddiThm}
	K_{\ell}(x,y)=\sum_{k=1}^{N(d,\ell)}Y_{\ell,k}(x)Y_{\ell,k}(y)=\frac{N(d,\ell)}{\omega_d}P_{\ell}(d+1;x\cdot y).
\end{equation*}


\section{Spherical convolution approximation}
\label{sec:sphQI}
Convolution operators are  powerful tools employed in  many  powerful approximation methods, as demonstrated in previous works such as \cite{smaleZhou2007-ConstrApprox-learning,wu2005generalized,ye_MAMS2019_generalized}. In this paper, we investigate a new type of spherical convolution operators $\ConvPhi$:
\begin{equation}\label{eq:ConvOper}
	(\ConvPhi f)(x)=(f*\varphi_{\rho})(x)=\int_{\mathbb{S}^d}f(y)\varphi_{\rho}(x\cdot y)\d\mu(y).
\end{equation}
Here the zonal kernel $\varphi_{\rho}$ is of the form as in \eqref{eq:DefiSphKernel}, which admits the following Fourier-Legendre expansion:
\begin{equation}\label{eq:expreScaleZonal}
	\scaleKer(x\cdot y)=\sum_{\ell=0}^{\infty}\widehat{\scaleKer}(\ell)K_{\ell}(x,y).
\end{equation}

The first task at hand is to identify conditions on the Fourier-Legendre coefficients  $\widehat{\scaleKer}(\ell)$
under which the convolution operators in \eqref{eq:ConvOper} deliver an ideal approximation power. 
For the two-dimensional sphere $\mathbb{S}^2$, Gomes, Kushpel and Levesley \cite{gomes2001approximation} investigated band-limited kernels under specific conditions employed in spherical convolution operators. We generalize this setting to accommodate scaled zonal kernels on spheres of integer dimensions $\mathbb{S}^d$ that satisfy the following assumption.
\begin{assumption}\label{Assump1:Kernel}
	There exist $s \ge 0$, $0<\rho_0<1$,  and  $C>0$ independent of $\ell$ and $\rho$ such that 
	\begin{equation}\label{eq:AssumpOnKernel}
		|1-\widehat{\scaleKer}(\ell)|\leq C\ell^s\rho^s, \quad 0<\rho \leq \rho_0, ~0\leq\ell\leq \Nro:=\Big\lfloor \frac{1}{\rho}-1\Big\rfloor,
	\end{equation}
    and
    \begin{equation}\label{eq:AssumpOnKernel1}
        |\widehat{\scaleKer}(\ell)|\leq C,  \quad \ell > \Nro.
    \end{equation}
\end{assumption}


\begin{remark} The special case $s=0$ amounts to requiring that the sequence
	$\{\widehat{\scaleKer}(\ell)\}_{\ell=0}^\infty$ be uniformly bounded with respect to $\rho$ and $\ell.$
	For the filtered hyperinterpolation kernel $\Phi_{L,G_a}$  defined in \eqref {hyperkernel}, we set $\rho = 1/L$ to get $ 	G_a(\rho \ell)=1$, for all $\ell \leq 1/\rho$ and $a\geq 1$.
	Thus,  filtered hyperinterpolation kernels satisfy Assumption \ref{Assump1:Kernel} with any $s>0$. Moreover, we will show in the sequel that  the spherical version of Poisson kernel satisfies Assumption \ref{Assump1:Kernel} with $s=1$ and the restriction to spheres of scaled Gaussian kernels and a certain type of compactly-supported kernels with $s=2$, respectively. We will also construct zonal kernels satisfying Assumption \ref{Assump1:Kernel} with higher orders.  This underscores that our spherical convolution framework extends the filtered hyperinterpolation method by relaxing kernel conditions. This extension accommodates a broader range of zonal kernels, which is advantageous in the context of multiresolution analysis, as the tuning parameter provides additional flexibility.
\end{remark}

\begin{theorem}\label{thm:SphConvConvergence}
	Suppose that the scaled zonal kernel $\varphi_{\rho}$ satisfies Assumption \ref{Assump1:Kernel}.
	Then for $0\leq \tau \leq s$ and $f\in H^{s}(\Sph^d)$, there is a constant $C$ independent of $\rho$ such that
	\begin{equation}\label{eq:ConvError}
		\|f-\ConvPhi f\|_{H^\tau}(\S^d)\leq C\rho^{s-\tau}\|f\|_{H^s(\mathbb{S}^d)}.
	\end{equation}
\end{theorem}	
\begin{proof}
	Using \eqref{eq:FunkHeckeFormula}
	and \eqref{eq:fexpansion}, we write 
	\begin{equation}\label{eq:ExpanOfConvOperator}
		(\ConvPhi f)(x)=\sum_{\ell=0}^{\infty}\sum_{k=1}^{N(d,\ell)}\widehat{f}_{\ell,k}\ConvPhi Y_{\ell,k}(x)=\sum_{\ell=0}^{\infty}\sum_{k=1}^{N(d,\ell)}\widehat{f}_{\ell,k}\widehat{\scaleKer}(\ell)Y_{\ell,k}(x).
	\end{equation}
	It then follows that
	\begin{equation}
		\begin{aligned}
			\|f-\ConvPhi f\|_{H^{\tau}(\Sph^d)}^2=&~\sum_{\ell=0}^{\infty}\sum_{k=1}^{N(d,\ell)}(1+\ell)^{2\tau}|\widehat{f}_{\ell,k}-\widehat{f}_{\ell,k}\widehat{\scaleKer}(\ell)|^2\\
			=&~\sum_{\ell=0}^{\infty}|1-\widehat{\scaleKer}(\ell)|^2\sum_{k=1}^{N(d,\ell)}(1+\ell)^{2\tau}|\widehat{f}_{\ell,k}|^2.
		\end{aligned}
	\end{equation}
	We split the infinite series on the right hand side in two parts: $\ell\leq N_{\rho}$ and $\ell>N_{\rho}$, and estimate them separately. We bound the first part using Assumption \ref{Assump1:Kernel} on $\varphi_{\rho}$,
	\begin{equation}\label{eq:ConvErrA}
		\begin{aligned}
			A~:=&~\sum_{\ell\leq \Nro}|1-\widehat{\scaleKer}(\ell)|^2\sum_{k=1}^{N(d,\ell)}(1+\ell)^{2\tau}|\widehat{f}_{\ell,k}|^2\\
			\leq&~C\sum_{\ell\leq \Nro}\ell^{2s}\rho^{2s}\sum_{k=1}^{N(d,\ell)}(1+\ell)^{2\tau}|\widehat{f}_{\ell,k}|^2\\
			=&~C\rho^{2s-2\tau}\sum_{\ell\leq \Nro}\ell^{2s}\rho^{2\tau}(1+\ell)^{2\tau-2s}\sum_{k=1}^{N(d,\ell)}(1+\ell)^{2s}|\widehat{f}_{\ell,k}|^2.
		\end{aligned}
	\end{equation}
	For $\ell\leq \Nro$, we have
	\begin{eqnarray*}
		\rho(\ell+1)\leq 1~~ \text{and}~~\ell^{2s}\rho^{2\tau}(1+\ell)^{2\tau-2s}\leq 1.
	\end{eqnarray*}
	We then plug this into \eqref{eq:ConvErrA} to get
	$$A\leq~C\rho^{2s-2\tau}\sum_{\ell\leq \Nro}\sum_{k=1}^{N(d,\ell)}(1+\ell)^{2s}|\widehat{f}_{\ell,k}|^2\leq C\rho^{2(s-\tau)}\|f\|_{H^{s}(\Sph^d)}^2.$$
	To bound the second part,
	we use Jensen's inequality to write
	\begin{equation*}
		\begin{aligned}
			B~:=&~\sum_{\ell> \Nro}|1-\widehat{\scaleKer}(\ell)|^2\sum_{k=1}^{N(d,\ell)}(1+\ell)^{2\tau}|\widehat{f}_{\ell,k}|^2\\
			\leq&~2\sum_{\ell> \Nro}\sum_{k=1}^{N(d,\ell)}(1+\ell)^{2\tau}|\widehat{f}_{\ell,k}|^2+2\sum_{\ell> \Nro}|\widehat{\scaleKer}(\ell)|^2\sum_{k=1}^{N(d,\ell)}(1+\ell)^{2\tau}|\widehat{f}_{\ell,k}|^2=:2(B_1+B_2).
		\end{aligned}
	\end{equation*}
	Since $\tau\leq s$, we have
	\begin{equation*}
		\begin{aligned}
			B_1\leq  \Nro^{2(\tau-s)}\sum_{\ell> \Nro}\sum_{k=1}^{N(d,\ell)}(1+\ell)^{2s}|\widehat{f}_{\ell,k}|^2
			\leq  C\rho^{2(s-\tau)}\|f\|_{H^{s}(\Sph^d)}^2.
		\end{aligned}
	\end{equation*}
	Using the uniform boundedness of sequence $\widehat{\scaleKer}(\ell)$, we derive
	\begin{equation*}
		\begin{aligned}
			B_2~=&~\sum_{\ell> \Nro}|\widehat{\scaleKer}(\ell)|^2\sum_{k=1}^{N(d,\ell)}(1+\ell)^{2\tau}|\widehat{f}_{\ell,k}|^2\\
			~\leq &~ C\rho^{2(s-\tau)}\sum_{\ell> \Nro}|\widehat{\scaleKer}(\ell)|^2\sum_{k=1}^{N(d,\ell)}(1+\ell)^{2s}|\widehat{f}_{\ell,k}|^2\\
			~\leq &~ C\rho^{2(s-\tau)}\|f\|_{H^{s}(\Sph^d)}^2,\\
		\end{aligned}
	\end{equation*}
	where in the last inequality we have used the assumption $|\widehat{\scaleKer}(\ell)|\leq C$ for $\ell> \Nro$. Combining all results completes the proof of the desired inequality \eqref{eq:ConvError}. 
\end{proof}

\subsection{Scaled radial kernels restricted on the sphere}
In this part, we investigate the restriction on spheres of  scaled radial kernels, including the spherical version of Poisson kernel (or simply Poisson kernel), Gaussian kernel, and compactly-supported radial kernels and derive closed formulas for their Fourier-Legendre coefficients \cite{mhaskar-1999Adv-approximation,narcowich2007approximation}. 

We will adopt the following transform pair (Fourier transform and its inverse) for a function $f\in L_2(\R^d)$: 
\begin{fequation}\label{eq:FourTranDef}
	\widehat{f}(\omega)=\int_{\R^d}f(x)e^{-\i \omega\cdot x}\d x, ~~\text{and}~~f(x)=\frac{1}{(2\pi)^d}\int_{\R^d}\widehat{f}(\omega)e^{\i \omega\cdot x}\d \omega.
\end{fequation}
Furthermore, we will assume that the Fourier transform and its inverse operator have been extended in a standard way so that they become applicable to the Schwartz class distributions, which we will refer to as  ``distributional Fourier transform". 
It is well-known that the Fourier transform of a radial function is still radial. Suppose $\phi\in L_1(\R^d)$ is radial, then its Fourier transform is represented by
\begin{fequation}\label{eq:RBFFourTran}
	\widehat{\phi}(\omega)=\mathcal{F}_d\phi(r)=(2\pi)^{d/2}r^{-(d-2)/2}\int_{0}^{\infty}\phi(t)t^{d/2}J_{(d-2)/2}(rt)\d t,~~r=\|\omega\|,
\end{fequation}
where $J_{\nu}(z)$ is the order-$\nu$ Bessel function of the first kind. Let $C^*(\R^{d+1})$ denote the totality of all continuous functions on $\R^{d+1}$ whose distributional Fourier transform $\widehat{\phi}$ is measurable on $\R^{d+1}$, and satisfies the following:
\begin{equation}\label{measure-condition}
	\int_{0\le |\xi|\le 1} |\xi|^{d-1}|\widehat{\phi}(\xi)|\d\xi<\infty \quad {\rm
		and}\quad \int_{|\xi|\ge 1} |\widehat{\phi}(\xi)|\d\xi<\infty.
\end{equation}
In  \cite{narcowich2007approximation}, the authors established a relation between the Fourier-Legendre coefficients of zonal functions and the Fourier transform of radial functions. 
\begin{lemma}{(\cite[Proposition 3.1]{narcowich2007approximation})}\label{lem:FourierCoeff}
	Let $\phi \in C^*(\R^{d+1}) $ be radial and $\varphi(x\cdot y):=\phi(x-y)|_{x,y\in\S^d}$. Then the following identities hold true. 
	\begin{equation*}
		\widehat{\varphi}(\ell)=\int_{0}^{\infty} t\widehat{\phi}(t)J_v^2(t)\mathrm{d} t,~~\nu:=\ell+\frac{d-1}{2}.
	\end{equation*}
\end{lemma}

We will use Lemma \ref{lem:FourierCoeff} to calculate the Fourier-Legendre coefficients of zonal kernels derived from scaled radial kernels.

\subsubsection{Poisson kernel}  Consider the Poisson kernel on the unit sphere $\mathbb{S}^{d}$, given by
\begin{equation}\label{eq:PoissonKernel}
	P_{\alpha}(x\cdot y)=\frac{1-\alpha^2}{\omega_d(1+\alpha^2-2\alpha x\cdot y)^{(d+1)/2}}, ~~0\leq \alpha<1,~ x,y\in\mathbb{S}^d,
\end{equation}
where $\omega_{d}$ is the surface area of $\S^d$. We view $P_{\alpha}$ as the restriction to $\S^d \times \S^d$ of the radial kernel $\phi(x,y;\alpha)$ on $\R^{d+1} \times \R^{d+1}$: 
$$\phi(x,y;\alpha)=\frac{1-\alpha^2}{\omega_d\big((1-\alpha)^2+\alpha\|x-y\|^2\big)^{(d+1)/2}},~~0\leq \alpha<1,$$
which reduces to 
\eqref{eq:PoissonKernel} when we enforce the conditions $\|x\|=\|y\|=1$. 
The Poisson kernel enjoys the following simple
expansion in terms of spherical harmonics, as detailed in M\"{u}ller \cite[Lem.~17]{muller1966SphHarm},
\begin{equation}\label{eq:CoeffPoisson}
	\varphi_{\alpha}^P(x\cdot y)=\sum_{\ell=0}^{\infty}\alpha^{\ell}\sum_{k=1}^{N(d,\ell)}Y_{\ell,k}(x)Y_{\ell,k}(y),~~0\leq\alpha<1.
\end{equation}
Upon the substitution $\alpha=1-\rho$ and an application of the inequality
\[
1-(1-\rho)^\ell \le \ell\ \rho, \quad 0<\rho<1, \quad \ell \in \mathbb N,
\]
we see that Poisson kernel satisfies Assumption \ref{Assump1:Kernel} with $s=1.$

\subsubsection{Gaussian kernel}

We define the Gaussian kernel $G_{\rho}(x):=\exp(-\frac{\|x\|^2}{2\rho^2})$, 
then its restriction onto the unit sphere $\S^d$ takes the form:
\begin{equation}\label{eq:GaussianSBF}
	G_\rho(x\cdot y)=\exp\Big(-\frac{2-2x\cdot y}{2\rho^2}\Big)=\exp\Big(-\frac{1-x\cdot y}{\rho^2}\Big).
\end{equation}
By defining the scaled zonal kernel as
\begin{equation}\label{eq:SphGauss}
	\varphi_{\rho}(x\cdot y):=(2\pi)^{-d/2}\rho^{-d}G_\rho(x\cdot y),
\end{equation}
we can obtain the following property.
\begin{proposition}\label{prop:AsympExapGaussian}
	Let $0<\rho<1$, we have the following expansion
	\begin{equation}\label{eq:HadmardSeries}
		\widehat{\varphi_{\rho}}(\ell)=\sum_{j=0}^{\infty}a_{j}(\nu)\rho^{2j}\cdot P(\nu+j+\frac{1}{2},\frac{2}{\rho^2}),~~\nu=\ell+(d-1)/2,
	\end{equation}
	where  $P(b,x)=\gamma(b,x)/\Gamma(b)$  denotes the normalized incomplete gamma function and the expansion coefficients are given by
    \begin{equation}\label{eq:ajv}
        a_{j}(\nu)=\frac{(\frac{1}{2}-\nu)_{j}(\frac{1}{2}+\nu)_{j}}{2^{j} j!}=\frac{(-1)^j(4\nu^2-1)(4\nu^2-9)\cdots[4\nu^2-(2j-1)^2]}{8^j\cdot j!}
    \end{equation}
	with $(x)_{j}=\Gamma(x+j)/\Gamma(x)$ denoting the Pochhammer symbol.  In particular, if $d$ is even, i.e., there is a $k\in\mathbb{N}$ such that $d=2k$,  the expansion \eqref{eq:HadmardSeries} reduces to the finite sum,
	$$\widehat{\varphi_{\rho}}(\ell)=\sum_{j=0}^{\ell+k-1}a_j(\nu)\rho^{2j}\cdot P(\nu+j+\frac{1}{2},\frac{2}{\rho^2}).$$

\end{proposition}
\begin{proof}
	The Fourier-Legendre coefficients of \eqref{eq:GaussianSBF} were originally derived in \cite{narcowich2007approximation},
	\begin{equation}\label{eq:FourCoeff-GSBF}
		\widehat{G_{\rho}}(\ell)=(2\pi)^{(d+1)/2}\ \rho^{d-1}\exp\Big(-\frac{1}{\rho^2}\Big)I_{\nu}\Big(\frac{1}{\rho^2}\Big),~~\nu=\ell+(d-1)/2,
	\end{equation}
	where $I_{\nu}$ represents the modified Bessel function of the first kind, which enjoys the following Hadamard expansion (see \cite[7.25,~(1)]{watson1966treatise}):
	\begin{equation*}
		I_{\nu}(z)= \frac{e^{z}}{\sqrt{2\pi z}}\sum_{j=0}^{\infty}\frac{(\frac{1}{2}-\nu)_j}{j!(2z)^j}\frac{\gamma(\nu+j+\frac{1}{2},2z)}{\Gamma(\nu+\frac{1}{2})},
	\end{equation*}
	where $\gamma$ denotes the incomplete gamma function. 
	We furthermore use the formula \cite[(2.5)]{paris2001usePRSLA} to rewrite the above expansion as
	\begin{equation}\label{paris}
		I_{\nu}(z)=\frac{e^{z}}{\sqrt{2\pi z}}\sum_{j=0}^{\infty}\frac{a_{j}(\nu)}{z^j}\cdot P(\nu+j+\frac{1}{2},2z),
	\end{equation}
	where 
	$$P(b,x) = \frac{\gamma(b,x)}{\Gamma(b)},~~\gamma(b,x)=\int_{0}^x e^{-t}t^{b-1}\d t,~~b>0.$$
	Upon the substitution $z=\frac{1}{\rho^2}$ in \eqref{paris}, then \eqref{eq:FourCoeff-GSBF} becomes
	\begin{equation*}
		\widehat{G_{\rho}}(\ell)=(2\pi)^{d/2}\rho^d\sum_{j=0}^{\infty}a_{j}(\nu)\rho^{2j}\cdot P(\nu+j+\frac{1}{2},\frac{2}{\rho^2}),
	\end{equation*}
	which leads to the desired formula \eqref{eq:HadmardSeries}.
	If $d=2k$ for some $k \in \mathbb N$,  we have 
	$$(\frac{1}{2}-\nu)_j=\big(\frac{1}{2}-\ell-\frac{d-1}{2}\big)_j=(1-\ell-k)_j.$$
	Thus $a_j(\nu)=0$ for $j\geq \ell+k$ as shown in the following equation:
	\begin{equation*}
		(-n)_j=\left\{\begin{aligned}
			&\frac{(-1)^j\cdot n!}{(n-j)!},~&0\leq j\leq n,\\
			&0,~& j\geq n+1.
		\end{aligned}\right.
	\end{equation*}
	This completes the proof.
\end{proof}

\begin{proposition}\label{prop:GaussAssump2}
	The scaled zonal kernel  $\varphi_{\rho}$ defined by \eqref{eq:SphGauss} satisfies \eqref{eq:AssumpOnKernel} with $s=2$ for $0\leq \ell\leq N_{\rho}$.
\end{proposition}

\begin{proof}
	The asymptotic expansion of the modified Bessel function $I_{\nu}(z)$ has been extensively investigated in the literature; see e.g. \cite{Abramowitz1964handbook,nemes2017AAMerror,watson1966treatise}. In a recent work \cite[(95)]{nemes2017AAMerror}, the authors gave the following comprehensive formula for large argument expansion of $I_{\nu}(z)$:
	\begin{equation}\label{eq:ExpanRemainder}
    \begin{aligned}
        I_{\nu}(z)=&~\frac{e^{z}}{\sqrt{2\pi z}}\Big[\sum_{j=0}^{N-1}\frac{a_j(\nu)}{z^j}+\mathcal{R}_N^{(K)}(ze^{\mp\pi\i},\nu)\Big]\\
        &~~\pm\i e^{\pm\pi\i \nu}\frac{e^{-z}}{\sqrt{2\pi z}}\Big[\sum_{k=0}^{M-1}\frac{(-1)^ka_k(\nu)}{z^k}+\mathcal{R}_M^{(K)}(z,\nu)\Big],
    \end{aligned}
	\end{equation}
	in which the definition of $a_j(\nu)$ differs from ours by a factor of  $(-1)^j$, and $\mathcal{R}_N^{(K)}(z,\nu)$ is the remainder term with the error bound 
	\begin{equation}\label{gaussian-remainder-bound}
		|\mathcal{R}_N^{(K)}(z,\nu)|\leq  C_N\frac{|a_N(\nu)|}{|z|^N},
	\end{equation}
	where the coefficient $C_N$ is explicitly given in \cite[(75)-(76)]{nemes2017AAMerror}. In particular, for $|\arg (z)| < \frac{\pi}{2}$, $C_N=1.$ 
	Based on the above expansion, we derive the following equation for $\widehat{\varphi_{\rho}}(\ell)$:
	\begin{equation}
    \begin{aligned}
        \widehat{\varphi_{\rho}}(\ell)=&~\sum_{j=0}^{N-1}a_j(\nu)\rho^{2j}+\mathcal{R}_N^{(K)}\big(-\frac{1}{\rho^2},\nu\big)\\
        &~\pm\i e^{\pm\pi\i \nu}e^{-\frac{2}{\rho^2}}\Big[\sum_{k=0}^{M-1}(-1)^ka_k(\nu)\rho^{2k}+\mathcal{R}_M^{(K)}(\frac{1}{\rho^2},\nu)\Big].
    \end{aligned}
	\end{equation}
	Setting $M=N$ in the above equation leads to the following error bound:
	\begin{equation}\label{eq:ErrBoundGaussian1}
		\begin{aligned}
			|\widehat{\varphi_{\rho}}(\ell)-1|
			\leq\Big|\sum_{j=0}^{N-1}a_j(\nu)\rho^{2j}-1\Big|+2|a_N(\nu)|\rho^{2N}+e^{-\frac{2}{\rho^2}}\Big|\sum_{k=0}^{N-1}(-1)^ka_k(\nu)\rho^{2k}\Big|,
		\end{aligned}
	\end{equation}
	where we have used the estimate from \cite[(75)]{nemes2017AAMerror} that $\mathcal{R}_N^{(K)}\big(\pm\frac{1}{\rho^2},\nu\big)\leq |a_N(\nu)|\rho^{2N}$.
	
	The ratio of two successive coefficients $a_j(\nu)$ is
	$$r_j(\nu):=\frac{a_{j+1}(\nu)}{a_{j}(\nu)}=-\frac{4\nu^2-(2j+1)^2}{8(j+1)},~j=0,1,2,\cdots.$$
	We observe that for $j > \nu - \frac{1}{2}$, the ratio $r_j(\nu)$ is positive, indicating that all the coefficients $a_j(\nu)$ have the same sign. However, for $0 \leq j \leq \nu - \frac{1}{2}$, the ratio $r_j(\nu)$ is negative, implying that the successive coefficients $a_j(\nu)$ have alternating signs.
	
	With this premise, we let $N=\lfloor \nu-\frac{1}{2}\rfloor+1$ in \eqref{eq:ErrBoundGaussian1} and denote 
	\begin{equation*}
		\begin{aligned}
			\mathcal{S}_{N}(\nu,\rho)=\sum_{j=0}^{N-1}(-1)^j\eta_j(\nu,\rho),~~\eta_j(\nu,\rho):=(-1)^ja_j(\nu)\rho^{2j},
		\end{aligned}
	\end{equation*}
	where it is evident that $\eta_j(\nu,\rho)>0$ since $a_j(\nu)$ has alternating signs. Furthermore, considering the ratio of successive terms, we find
	\begin{equation*}
		\begin{aligned}
			\frac{\eta_{j+1}(\nu,\rho)}{\eta_{j}(\nu,\rho)}=&~\frac{4\nu^2-(2j+1)^2}{8(j+1)}\rho^2
			<\frac{\nu^2\rho^2}{2(j+1)}<1.
		\end{aligned}
	\end{equation*}
	The final inequality hinges on the assumption that $\nu\rho < 1$. This shows that $\eta_j$ is positive and decreases monotonically with respect to $j$, which implies that \begin{equation}\label{eq:estimateS}
		\begin{aligned}
			\big|\mathcal{S}_{N}(\nu,\rho)-\eta_0(\nu,\rho)\big|\leq&~  |\eta_1(\nu,\rho)|.
		\end{aligned}
	\end{equation} 
	That is
	\begin{equation}\label{eq:CoeffEst1}
		\big|\mathcal{S}_{N}(\nu,\rho)-1\big|< \frac{1}{2}\nu^2\rho^2.
	\end{equation}
	We further obtain the following estimates in a straight-forward fashion:
	\begin{equation}\label{eq:CoeffEst2}
		|a_N(\nu)|=\Big|\frac{(4\nu^2-1)(4\nu^2-9)\cdots[4\nu^2-(2N-1)^2]}{8^N\cdot N!}\big|< \frac{\nu^{2N}}{2^N\cdot N!},
	\end{equation}
	and
	\begin{equation}\label{eq:CoeffEst3}
		\Big|\sum_{k=0}^{N-1}(-1)^ka_k(\nu)\rho^{2k}\Big|\leq \sum_{k=0}^{N-1}|a_k(\nu)|\rho^{2k}<\sum_{k=0}^{N-1}\frac{\nu^{2k}\rho^{2k}}{2^k\cdot k!}<\sum_{k=0}^{N-1}\frac{1}{2^k\cdot k!}<2.
	\end{equation}
	Incorporating \eqref{eq:CoeffEst1}-\eqref{eq:CoeffEst3} into \eqref{eq:ErrBoundGaussian1}, we derive
	$$|\widehat{\varphi_{\rho}}(\ell)-1|\leq \nu^2\rho^2+2e^{-\frac{2}{\rho^2}}.$$
	Finally, we bound the last term on the right hand side of the above inequality by $\rho^2$ to complete the proof.
\end{proof}

\subsubsection{Compactly-supported radial kernels} In this part, we are concerned with a class of compactly-supported radial kernels. Let $m>-1$, $\rho>0$, set the radial kernel $\phi_{\rho}^m(x)$ in $\R^{d+1}$ as
\begin{fequation}\label{eq:CSkernel}
	\phi_{\rho}^m(x)=\left\{
	\begin{aligned}
		&(1-\|x\|^2/\rho^2)^m,~~&\text{if}~\|x\|\leq \rho,\\
		&0,~~&\text{if}~\|x\|> \rho.
	\end{aligned}\right.
\end{fequation}
Let $\varphi_{\rho}^m(x\cdot y)$ be the associated scaled zonal kernel of $\phi_{\rho}^m$. We define a modified kernel 
\begin{equation}\label{eq:ScaledCompactKernel}
	\varphi_{\rho}(x\cdot y):=\frac{\Gamma(m+\frac{d}{2}+1)}{\pi^{\frac{d}{2}}\Gamma(m+1)\rho^d}\varphi_{\rho}^m(x\cdot y),~~x,y\in\S^d.
\end{equation}
 The Fourier-Legendre coefficients of the kernel are computed in the following proposition. 
\begin{proposition}\label{prop:CSkernel1}
	Let the scaled zonal kernel  $\varphi_{\rho}$ be defined by \eqref{eq:ScaledCompactKernel}, then we have
	\begin{equation}\label{eq:HyperExpanOfCSkernel}
		\widehat{\varphi_{\rho}}(\ell)=\hgeom{2}{1}(\ell+\frac{d}{2},-\ell-\frac{d}{2}+1;m+\frac{d}{2}+1;\frac{\rho^2}{4})=\sum_{j=0}^{\infty}\frac{a_j(\nu)}{(m+\frac{d}{2}+1)_j\cdot 2^j}\rho^{2j},
	\end{equation}
	where $\hgeom{2}{1}(a,b;c;z)$ is the hypergeometric function and $a_j(\nu)$ is defined in \eqref{eq:ajv}. If $d=2k$ ($k\in\mathbb{N}$), then the above infinite series reduces to the following finite sum:
	\begin{equation}\label{eq:evenDimensioanCSkernel}
		\widehat{\varphi_{\rho}}(\ell)~=\sum_{j=0}^{\ell+k-1}\frac{a_j(\nu)}{(m+k+1)_j \cdot 2^j}\rho^{2j},~~\nu=\ell+\frac{d-1}{2}.
	\end{equation}
\end{proposition}

To facilitate the calculation, we will outline some key properties of the hypergeometric function, which will be essential for the subsequent analysis.
\begin{definition}\label{prop:HypergeomFun}
The hypergeometric function is defined as
\begin{equation}
\begin{aligned}
	\hgeom{2}{1}(a,b;c;z)=&\sum_{k=0}^{\infty}\frac{(a)_k(b)_k}{(c)_kk!}z^k=\frac{\Gamma(c)}{\Gamma(a)\Gamma(b)}\sum_{k=0}^{\infty}\frac{\Gamma(a+k)\Gamma(b+k)}{\Gamma(c+k)k!}z^k\\
	=&1+\frac{ab}{c}z+\frac{a(a+1)b(b+1)}{c(c+1)2!}z^2+\cdots.
\end{aligned}
\end{equation}
If $a=-n$, $n \in \mathbb{N}$ and $c\neq0,-1,-2,\ldots$, $\hgeom{2}{1}(a,b;c;z)$ is a polynomial,
\begin{equation}\label{eq:negativeNumberHypergeoFunc}
\hgeom{2}{1}(-n,b;c;z)=\sum_{k=0}^n\frac{(-n)_k(b)_k}{(c)_k k!}z^k=\sum_{k=0}^{n}(-1)^k{n\choose k}\frac{(b)_k}{(c)_k }z^k.
\end{equation}
The Euler's hypergeometric transformation is given by
\begin{equation}\label{eq:EulerTransformation}
\hgeom{2}{1}(a,b;c;z)=(1-z)^{1-a-b}\hgeom{2}{1}(c-a,c-b;c;z).
\end{equation}
\end{definition}

To prove Proposition \ref{prop:CSkernel1}, we begin by presenting some useful lemmas. The Fourier transform of compactly supported kernels is given as follows \cite[Thm.~5.26]{wendland2004scattered}.
\begin{lemma}\label{lem:FourTranCS}
The Fourier transform of $\phi_{\rho}^{m}$ as a function in $\R^{d+1}$ is given by
\begin{equation*}
\widehat{\phi}_{\rho}^{m}(r)=(2\pi)^{(d+1)/2}2^{m}\Gamma(m+1)r^{-(d+1)/2-m}\rho^{(d+1)/2-m}J_{(d+1)/2+m}(\rho r).
\end{equation*}
\end{lemma}

Next, we establish a useful integral formula involving triple products of Bessel functions.
\begin{lemma}\label{lem:IntegralThreeBessel}
For $\text{Re}(\mu)>-\frac{1}{2}$, $\text{Re}(\nu)>-\frac{1}{2}$ and $a$, $b$, $c$ representing the three sides of a triangle, let $x=\frac{b^2+c^2-a^2}{2bc}$, it holds that
\begin{equation}\label{eq:Int3Bessel}
\begin{aligned}
    \int_{0}^{\infty}J_{\mu}(at)J_{\nu}(bt)J_{\nu}(ct)t^{1-\mu}dt
=&~\frac{(bc)^{\mu-1}}{(2a)^{\mu}\Gamma(\mu+\frac{1}{2})\Gamma(\frac{1}{2})}(1-x^2)^{\mu-\frac{1}{2}}\\
&~~\cdot\hgeom{2}{1}(\mu-\nu,\mu+\nu;\mu+\frac{1}{2};\frac{1-x}{2}).
\end{aligned}
\end{equation}
\end{lemma}

\begin{proof}
From \cite[Sec.~13.46]{watson1966treatise}, we obtain the following integral formula
\begin{equation}\label{eq:IntegralFormula1}
\begin{aligned}
    \int_{0}^{\infty}J_{\mu}(at)J_{\nu}(bt)J_{\nu}(ct)t^{1-\mu}dt
=&~\frac{(2bc)^{\mu-1}}{a^{\mu}\Gamma(\mu+\frac{1}{2})\Gamma(\frac{1}{2})}\sin^{2\mu-1}\frac{A}{2}\\
&~~\cdot\hgeom{2}{1}(\frac{1}{2}+\nu,\frac{1}{2}-\nu;\mu+\frac{1}{2};\sin^2\frac{A}{2}),
\end{aligned}
\end{equation}
where $A=\arccos\frac{b^2+c^2-a^2}{2bc}$.
By making the substitution
\begin{equation*}
x=\cos A,~\sin^2\frac{A}{2}=\frac{1-\cos A}{2}=\frac{1-x}{2},
\end{equation*}
we can rewrite the expression as
\begin{fequation}\label{eq:IntegralFormula2}
\begin{aligned}
	&\sin^{2\mu-1}\frac{A}{2}\cdot\hgeom{2}{1}(\frac{1}{2}+\nu,\frac{1}{2}-\nu;\mu+\frac{1}{2};\sin^2\frac{A}{2})\\
	=&~\Big(\frac{1-x}{2}\Big)^{\mu-\frac{1}{2}}\cdot\hgeom{2}{1}(\frac{1}{2}+\nu,\frac{1}{2}-\nu;\mu+\frac{1}{2};\frac{1-x}{2})\\
	=&~\Big(\frac{1-x^2}{4}\Big)^{\mu-\frac{1}{2}}\cdot\hgeom{2}{1}(\mu-\nu,\mu+\nu;\mu+\frac{1}{2};\frac{1-x}{2}),
\end{aligned}
\end{fequation}
where the second equality follows from Euler's hypergeometric transformation, as stated in Proposition \ref{prop:HypergeomFun}. Substituting \eqref{eq:IntegralFormula2} into \eqref{eq:IntegralFormula1} yields the desired integral formula.
\end{proof}

\begin{proof}[\textbf{Proof of Proposition \ref{prop:CSkernel1}}]
We begin by applying Lemma \ref{lem:FourierCoeff} to the kernel $\phi_{\rho}^{m}$, which yields the following expression for the Fourier-Legendre coefficients
\begin{equation}\label{eq:Four_Tran_CS_ker}
\begin{aligned}
    \widehat{\varphi}_{\rho}^m(\ell)=&~(2\pi)^{(d+1)/2}2^{m}\Gamma(m+1)\rho^{(d+1)/2-m}\\
&~~\cdot\int_{0}^{\infty} t^{1-(d+1)/2-m}J_{(d+1)/2+m}(\rho t)J_{\nu}^2(t)dt.
\end{aligned}
\end{equation}
Next, we use Lemma \ref{lem:IntegralThreeBessel} with the parameters $\mu=m+\frac{d+1}{2}$, $\nu=\ell+\frac{d-1}{2}$, $a= \rho$, $b=c=1$. This allows us to express the integral as
\begin{equation}\label{eq:IntegralRepre2}
\begin{aligned}
	&\int_{0}^{\infty} t^{1-(d+1)/2-m}J_{(d+1)/2+m}(\rho t)J_{\nu}^2(t)dt\\
	=&~\frac{1}{(2\rho)^{m+\frac{d+1}{2}}\Gamma(m+\frac{d}{2}+1)\Gamma(\frac{1}{2})}(1-x^2)^{m+\frac{d}{2}}\\
	&~~\cdot\hgeom{2}{1}(m+1-\ell,m+\ell+d;m+\frac{d}{2}+1;\frac{1-x}{2}).
\end{aligned}
\end{equation}
In this context, as per Lemma \ref{lem:IntegralThreeBessel}, we consider the triple $(1,1,\rho)$ as the three sides of a triangle and easily derive the following equations: 
\begin{equation*}
x=\cos A=1-\frac{\rho^2}{2},~~~ \frac{1-x}{2}=\frac{\rho^2}{4},~~~1-x^2=\rho^2(1-\frac{\rho^2}{4}).
\end{equation*}
Substituting these into the integral representation, we obtain
\begin{equation}
\begin{aligned}
	&\int_{0}^{\infty} t^{1-(d+1)/2-m}J_{(d+1)/2+m}(\rho t)J_{\nu}^2(t)dt\\
	=&~\frac{2^{-m-\frac{d+1}{2}}\rho^{m+\frac{d-1}{2}}}{\Gamma(m+\frac{d}{2}+1)\Gamma(\frac{1}{2})}\Big(1-\frac{\rho^2}{4}\Big)^{m+\frac{d}{2}}\\
	&~~\cdot\hgeom{2}{1}(m+1-\ell,m+\ell+d;m+\frac{d}{2}+1;\frac{\rho^2}{4}).
\end{aligned}
\end{equation}
Therefore, combining above results back into \eqref{eq:Four_Tran_CS_ker}, we derive
\begin{equation}
\begin{aligned}
	\widehat{\varphi}_{\rho}^m(\ell)
	=&~\frac{\pi^{\frac{d}{2}}\Gamma(m+1)\rho^d}{\Gamma(m+\frac{d}{2}+1)}(1-\frac{\rho^2}{4})^{m+\frac{d}{2}}\\
    &~~\cdot\hgeom{2}{1}(m+1-\ell,m+\ell+d;m+\frac{d}{2}+1;\frac{\rho^2}{4}).
\end{aligned}
\end{equation}
By applying Euler's hypergeometric transformation \eqref{eq:EulerTransformation}, we simplify this to
\begin{fequation}\label{eq:FourCoeffCS2}
\widehat{\varphi}_{\rho}^m(\ell)=\frac{\pi^{\frac{d}{2}}\Gamma(m+1)\rho^d}{\Gamma(m+\frac{d}{2}+1)}\hgeom{2}{1}(\ell+\frac{d}{2},-\ell-\frac{d}{2}+1;m+\frac{d}{2}+1;\frac{\rho^2}{4}).
\end{fequation} 
Letting $\nu=\ell+\frac{d-1}{2}$, we can rewrite the hypergeometric function as
\begin{align*}
&\hgeom{2}{1}(\ell+\frac{d}{2},-\ell-\frac{d}{2}+1;m+\frac{d}{2}+1;\frac{\rho^2}{4})\\
=&~\hgeom{2}{1}(\nu+\frac{1}{2},\frac{1}{2}-\nu;m+\frac{d}{2}+1;\frac{\rho^2}{4})\\
=&~\sum_{j=0}^{\infty}\frac{(\frac{1}{2}-\nu)_j(\frac{1}{2}+\nu)_j}{(m+\frac{d}{2}+1)_j\cdot j!}\Big(\frac{\rho^2}{4}\Big)^j\\
=&~\sum_{j=0}^{\infty}\frac{a_j(\nu)}{(m+\frac{d}{2}+1)_j\cdot 2^j}\rho^{2j},
\end{align*}
which completes the proof.
\end{proof}

Finally, we can prove the following property.

\begin{proposition}\label{prop:CSAssump1}
	The scaled zonal kernel  $\varphi_{\rho}$ defined by \eqref{eq:ScaledCompactKernel} satisfies \eqref{eq:AssumpOnKernel} with $s=2$ for $0\leq \ell\leq N_{\rho}$.
\end{proposition}

\begin{proof}
    The expansion of hypergeometric function had been proved in \cite[Thm.~2.2]{nemes-2020SIAMMA-large},
\begin{equation}\label{eq:Hypergeom FuncExpan}
	\hgeom{2}{1}(a,b;c;z)=\sum_{j=0}^{N-1}\frac{(a)_j(b)_j}{(c)_j\cdot j!}z^n+R_N^{(F)}(z,a,b,c),
\end{equation}
with the following upper bound for the remainder:
\begin{equation}\label{eq:RemainderEstimate}
	|R_N^{(F)}(z,a,b,c)|\leq A_N\Big|\frac{(a)_N(b)_N}{(c)_N\cdot N!}z^N\Big|,
\end{equation}
where $A_N$ is explicitly given by 
$$A_N=\Big|\frac{N}{a}\Big|+\Big|\frac{\Gamma(c+N)(1+\frac{N}{a})K^a}{\Gamma(c-b)(c-b)^{N+b}}\Big|,~~K=\min\big(2,\max(1,\frac{1}{1-2z})\big),$$
provided that $a$, $b$, and $c$ are real numbers and $N>-b$.  Here, we have restricted the complex variable $z$ on the nonnegative part of the real axis.

By taking $N=\lfloor -b\rfloor+1$, we have $-b<N\leq-b+1$, $N\leq a\leq N+1$, and we can estimate that 
\begin{align*}
A_N\leq1+\Big|\frac{\Gamma(c+N)2^{N+2}}{\Gamma(c-b)(c-b)^{N+b}}\Big|
< 1+2^{N+2}\cdot (c+N-1).
\end{align*}
In the case of \eqref{eq:HyperExpanOfCSkernel}, we have $$a=\ell+\frac{d}{2}=\frac{1}{2}+\nu,~~b=\frac{1}{2}-\nu,~~c=m+\frac{d}{2}+1,~~z=\frac{\rho^2}{4}.$$ 
Substituting these values into $R_N^{(F)}(z,a,b,c)$, we obtain
\begin{align*}
|R_N^{(F)}(z,a,b,c)|\leq &~A_N\Big|\frac{(\frac{1}{2}+\nu)_N\cdot(\frac{1}{2}-\nu)_N}{(c)_N\cdot N!}\Big(\frac{\rho^2}{4}\Big)^N\Big|
= ~A_N\frac{a_N(\nu)}{(c)_N\cdot 2^N}\rho^{2N}\\
< &~\Big[1+2^{N+2}\cdot (c+N-1)\Big]\frac{\nu^{2N}\rho^{2N}}{(c)_N\cdot 2^{2N}\cdot N!}\\
=&~\Big[1+2^{N+2}\cdot (c+N-1)\Big]\frac{\nu^{2N}\rho^{2N}}{c\cdot(c+1)\cdot (c+N-1)\cdot 2^{2N}\cdot N!}\\
< &~\frac{\nu^{2N}\rho^{2N}}{2^{N-2}\cdot N!},
\end{align*}
where we have used the estimate of $a_N(\nu)$ in \eqref{eq:CoeffEst2}. Thus we can use a similar method as in Proposition \ref{prop:GaussAssump2} to complete the proof.
\end{proof}

\subsection{High-order kernels} 
\label{subsec:HighOrder}

To construct zonal kernels that satisfy Assumption \ref{Assump1:Kernel} with large $s$, we employ strategies analogous to those used in the construction of radial kernels in Euclidean spaces. Such strategies often involve taking derivatives or forming linear combinations of scaled radial kernels to meet high-order generalized Strang-Fix conditions, as discussed in various studies \cite{beale-Majda1985high,mazya-schimidt2001quasi,gao2017constructing,franz-wendland2018convergence,ramming-wendland2018kernel,gao-zhou2020multiscale,franz-wendland2023multilevel}, which are also referred to as moment conditions in some contexts; see e.g. \cite{franz-wendland2023multilevel}. In the current zonal-kernel setting, we find that the ``linear combination" approach is more effective. To illustrate,  we will choose the compactly supported kernel \eqref{eq:ScaledCompactKernel} as a model, even though the same procedure works for many other zonal kernels.

\begin{proposition}
	Fix a $K \in \mathbb N$ and $0<a_1<a_2<\cdots<a_K \leq 1$.      
	Let $\varphi_{\rho}$ be defined as in \eqref{eq:ScaledCompactKernel}. Set
	\begin{equation}\label{eq:highOrderSBF}
		\psi_{\rho}(x\cdot y)=\sum_{i=1}^K\lambda_i\varphi_{\rho_i}\big(x\cdot y\big), ~~ x,y \in \S^{d}, ~~\rho_i = a_i\rho,
	\end{equation}
	in which
	\begin{equation}\label{eq:HighOrderCoeff}
		\lambda_i=\prod_{j\neq i}^{K}\frac{a_j^2}{a_j^2-a_i^2},~i=1,\ldots,K.
	\end{equation} 
	Then for even dimension $d$, the zonal kernel $\psi_\rho$ satisfies Assumption \ref{Assump1:Kernel} with $s=2K$ for $0\leq \ell\leq N_{\rho}$.
\end{proposition}
\begin{proof}
	Referring to \eqref{eq:evenDimensioanCSkernel} (or Proposition \ref{prop:CSkernel1}), for even dimension ($d=2k$), we have 
	$$\widehat{\varphi_{\rho}}(\ell)=\sum_{j=0}^{\ell+k-1}b_j(\nu)\rho^{2j},~~b_j(\nu):=\frac{a_j(\nu)}{(m+k+1)_j \cdot 2^j}.$$
This together with  \eqref{eq:highOrderSBF} leads to 
	\begin{equation}\label{eq:CombExp-CSkernel}
		\begin{aligned}
			\widehat{\psi_{\rho}}(\ell)&~=\sum_{i=1}^K\lambda_i\widehat{\varphi_{\rho_i}}\big(\ell\big)\\
            &~=\sum_{i=1}^K\lambda_i b_0(\nu)+\sum_{i=1}^K\lambda_i\sum_{j=1}^{\ell+k-1} b_j(\nu)(a_i\rho)^{2j}\\
			&=~\sum_{i=1}^K\lambda_i+\sum_{j=1}^{\ell+k-1}b_j(\nu)\rho^{2j}\Big(\sum_{i=1}^K \lambda_ia_i^{2j}\Big)\\
          &   =~\sum_{i=1}^K\lambda_i+\sum_{j=1}^{\ell+k-1}\Big(\sum_{i=1}^K \lambda_ia_i^{2j}\Big)b_j(\nu)\rho^{2j}.
		\end{aligned}
	\end{equation}
	To show that $\widehat{\psi_{\rho}}$ satisfies Assumption \ref{Assump1:Kernel} with $s=2K$, we enforce the following condition:
	\begin{equation}    \label{sys-lambda}
		\sum_{i=1}^K\lambda_i=1;~~\sum_{i=1}^K\lambda_ia_i^{2j}=0,~j=1,\ldots,K-1.
	\end{equation}
	The above is  a  $K \times K$ linear system  of equations with unknowns $\lambda_1,\ldots,\lambda_K$ which enjoys 
     a Vandermonde coefficient matrix. By Cramer's rule, the system has a unique solution given explicitly as in \eqref{eq:HighOrderCoeff} under the assumption that the positive numbers $a_i$ are distinct.       
	Hence, \eqref{eq:CombExp-CSkernel} reduces to 
	\begin{equation}\label{eq:CombExp-CSkernel2}
		\widehat{\psi_{\rho}}(\ell)=\left\{
		\begin{aligned}
			&1,~~&\ell+k-1<K,\\
			&1+\sum_{j=K}^{\ell+k-1}b_j(\nu)\rho^{2j}\sum_{i=1}^K\lambda_ia_i^{2j},~~&\ell+k-1\geq K.
		\end{aligned}\right.
	\end{equation}
	With the notation  $A:={\displaystyle \max_{1\le i\leq K} a_i}\leq 1$, we have the following estimate:
	\begin{equation}\label{eq:CSKernelEstimate_EvenDim}
		\begin{aligned}
			|\widehat{\psi_{\rho}}(\ell)-1|\leq&~ A^{2(\ell-k+1)}\Big(\sum_{i=1}^K{|\lambda_i|}\Big)\cdot\sum_{j=K}^{\ell+k-1}\big|b_j(\nu)\big|\rho^{2j}\\
			\leq &~\Big(\sum_{i=1}^K{|\lambda_i|}\Big)\cdot|b_K(\nu)|\rho^{2K}\sum_{j=K}^{\ell+k-1}\Big|\frac{b_j(\nu)}{b_K(\nu)}\Big|\rho^{2j-2K}.
		\end{aligned}
	\end{equation}
	We estimate the ratio: 
	\begin{equation*}
		\Big|\frac{b_j(\nu)}{b_K(\nu)}\Big|= \frac{a_j(\nu)\cdot (m+k+1)_K}{a_K(\nu)\cdot(m+k+1)_j \cdot 2^{j-K}}< \frac{\nu^{2(j-K)}K!}{2^{j-K}j!}\frac{1}{2^{j-K}},~j\geq K.
	\end{equation*}
	It then follows that 
	\begin{equation*}
		\sum_{j=K}^{\ell+k-1}\Big|\frac{b_j(\nu)}{b_K(\nu)}\Big|\rho^{2j-2K}< \sum_{j=K}^{\ell+k-1}\frac{\nu^{2(j-K)}K!}{2^{j-K}j!}\frac{1}{2^{j-K}}\rho^{2j-2K}< \sum_{j=0}^{\ell+k-1-K}(\frac{1}{4}\nu\rho)^j<2,
	\end{equation*}
	where we have used $\nu\rho<1$. Substituting the above inequality into \eqref{eq:CSKernelEstimate_EvenDim} completes the proof. 
\end{proof}

\section{Spherical  quasi-interpolation}
\label{sec:sphqi_scheme}
Let $X:=\{x_j\}^{|X|}_{j=1}$ be a finite subset of $\S^d,$ where $|X|$ denotes the cardinality of $X.$ For $n \in \mathbb N$, we call 
\begin{equation*}
	\QuadSet:=\{(w_j,x_j)|w_j\in\R,~x_j\in\S^d,~j=1,2,\ldots |X|\},
\end{equation*}
a quadrature rule of order $n$ on $\S^d$ if 
\begin{equation}\label{eq:quadrature_formula}
	\IntSph p(x)\d\mu(x)=\sum_{j=1}^{|X|}w_j p(x_j), \quad  p \in\Pi_n(\S^d).
\end{equation}
The numbers $w_j$ are called weights. If  $w_j>0$  for all $j=1,\ldots,|X|$,  then we call $\QuadSet$ a positive quadrature rule, which has been studied extensively in the literature; see, e.g., \cite{hesse_2006JAT_cubature,hesse2010numerical-book,montufar_2022FoCom_distributed}. 
By applying a positive quadrature rule $\QuadSet$ to the convolution integral $\ConvPhi f=f*\varphi_{\rho}$ in  \eqref{eq:ConvOper}, we obtain a quasi-interpolant of $f$ in the form:
\begin{equation}\label{eq:Spherical_quasi_interp}
	\ConvPhi^Xf(x):=\sum_{j=1}^{|X|}w_jf(x_j)\varphi_{\rho}(x_j\cdot x).
\end{equation}

The approximation error of the spherical quasi-interpolant in \eqref{eq:Spherical_quasi_interp} can be decomposed into two components: the convolution error and the discretization error. The analysis of the convolution error is presented in Theorem \ref{thm:SphConvConvergence}. To estimate the discretization error, we make the following assumption on the scaled zonal kernel we employ for quasi-interpolation.

\begin{assumption}\label{Assump2:Kernel}
	There exist two positive constants $s$ and $C$, such that 
	$$\widehat{\scaleKer}(\ell)\leq C(1+\rho\ell)^{-2s}, \quad \forall\ \ell>N_{\rho},$$
	where $N_{\rho}$ is defined as in Assumption \ref{Assump1:Kernel}.
\end{assumption}


Literature abounds with studies of decay properties of the Fourier-Legendre coefficients of scaled zonal kernels; see \cite{narcowich2007approximation,gia_2010SINUM_multiscale,hubbert2002radial}. Specifically, Le Gia, Sloan and Wendland \cite{gia_2010SINUM_multiscale} showed that if a radial function has compact support and its Fourier transform satisfies $\widehat{\phi}(\omega)\sim(1+\|\omega\|)^{-2s-1}$ with $s>d/2$, then the associated scaled zonal kernel $\varphi_{\rho}$ has Fourier-Legendre coefficients satisfying $\widehat{\scaleKer}(\ell)\sim(1+\rho\ell)^{-2s}$. 

Suppose that the kernel $\scaleKer$ satisfies Assumption \ref{Assump2:Kernel} with $s>d/2$.  For $f\in H^{s}(\S^d)$, we define:
$$\|f\|_{\varphi_{\rho}}^2:=\sum_{\ell=0}^{\infty}\sum_{k=1}^{N(d,\ell)}\frac{|\widehat{f}_{\ell,k}|^2}{\widehat{\scaleKer}(\ell)}.$$
This conforms to the conventional definition of the native space induced by $\varphi_{\rho}$:  
$$\mathcal{N}_{\varphi_{\rho}}:=\{f\in \mathcal{D}'(\S^d):\|f\|_{\varphi_{\rho}}< \infty\},$$
where $\mathcal{D}'(\S^d)$ represents all tempered distributions defined on $\S^d$.  

\begin{theorem}\label{thm:IntegralError}
	Let $f\in H^{s}(\S^d)$, $s>d/2$ and 
	$\scaleKer\in L_1(\S^d)$ satisfy Assumptions \ref{Assump1:Kernel} and \ref{Assump2:Kernel}. Suppose that $X$ is a quasi-uniform set and that  $\QuadSet$ is a positive quadrature rule of order $2n$ with $n>N_{\rho}$.  Then there exists a constant $C$ depends only on $d$ and $s$ such that
	$$\|f*\varphi_{\rho}-\ConvPhi^Xf\|_{L_2(\S^d)}\leq C\big(n^{-s}+n^{-s+d/2}+n^{-s}\rho^{-s}\big)\|f\|_{H^{s}(\S^d)}.$$
\end{theorem}

The proof of Theorem \ref{thm:IntegralError} needs several lemmas
that may already be accessible in the existing literature. We included proofs  here for completeness.

\begin{lemma}\label{lem:L2Proj}
	Let $f\in H^{s}(\S^d)$ with $s>d/2$, and $\scaleKer$ be defined as \eqref{eq:DefiSphKernel} satisfying Assumption \ref{Assump2:Kernel}. Define $\mathcal{P}f$, $\P\scaleKer$ be the $L_2$ projections of $f$, $\scaleKer$ onto the spherical polynomial space $\Pi_n(\S^d)$ of degree $n>N_{\rho}$, respectively. Then we have the following estimates  
	$$\|\P f\|_{L_2(\S^d)}\leq \|f\|_{L_2(\S^d)},~~\|f-\P f\|_{L_2(\S^d)}\leq n^{-s}\|f\|_{H^{s}(\S^d)},$$
	and
	$$\|\scaleKer-\P \scaleKer\|_{L_2(\S^d)}\leq C n^{-s}\rho^{-s}.$$
\end{lemma}
\begin{proof}
	The $L_2$ projection gives the following equations and inequalities: $$\widehat{f}_{\ell,k}=\widehat{\P f}_{\ell,k},~~\ell\leq n,~~ k=1,\ldots, N(d,\ell),$$
	and $\|\P f\|_{L_2(\S^d)}\leq \|f\|_{L_2(\S^d)}$. It then follows that
	\begin{align*}
		\|f-\P f\|^2_{L_2(\S^d)}=&\sum_{\ell=0}^{\infty}\sum_{k=1}^{N(d,\ell)}|\widehat{f}_{\ell,k}-\widehat{\P f}_{\ell,k}|^2=\sum_{\ell>n}^{\infty}\sum_{k=1}^{N(d,\ell)}|\widehat{f}_{\ell,k}|^2\\
		\leq &~n^{-2s}\sum_{\ell>n}^{\infty}\sum_{k=1}^{N(d,\ell)}(1+\ell)^{2s}|\widehat{f}_{\ell,k}|^2
		\leq  n^{-2s}\|f\|^2_{H^{s}(\S^d)}.
	\end{align*}
	Similarly, we get
	\begin{align*}
		\|\scaleKer-\P \scaleKer\|^2_{L_2}
		=&~ \sum_{\ell>n}^{\infty}\sum_{k=1}^{N(d,\ell)}|\widehat{\scaleKer}(\ell)|^2
		\leq  n^{-2s}\sum_{\ell>n}^{\infty}\sum_{k=1}^{N(d,\ell)}(1+\ell)^{2s}|\widehat{\scaleKer}(\ell)|^2\\
		\leq &~n^{-2s}\sum_{\ell>n}^{\infty}\sum_{k=1}^{N(d,\ell)}\rho^{-2s}(1+\rho\ell)^{2s}|\widehat{\scaleKer}(\ell)|^2\\
		\leq&~   C n^{-2s}\rho^{-2s}\|\scaleKer\|_{\mathcal{N}_{\varphi_{\rho}}}^2\leq Cn^{-2s}\rho^{-2s}.
	\end{align*}
	This completes the proof.
\end{proof}

\begin{lemma}\label{lem:ProjPointWiseErr}
	Under the assumptions in Lemma \ref{lem:L2Proj}, the following inequality holds 
	\begin{align*}
		\|f-\P f\|_{L_{\infty}(\S^d)}
		\leq  C\omega_d^{-1/2}n^{-s+d/2}\|f\|_{H^{
				s}(\S^d)}.
	\end{align*}
\end{lemma}
\begin{proof} 
	Using \eqref{eq:fexpansion} we have
	\begin{equation}\label{eq:Differ_f_Pf}
		f(x_j)-\P f(x_j)=\sum_{\ell>n}^{\infty}\sum_{k=1}^{N(d,\ell)}\widehat{f}_{\ell,k}Y_{\ell,k}(x_j).
	\end{equation}
	By the addition theorem for spherical harmonics \eqref{eq:AddiThm}, we get
	$$\sum_{k=1}^{N(d,\ell)}|Y_{\ell,k}(x_j)|^2=\frac{N(d,\ell)}{\omega_d}.$$
	Hence, applying Cauchy-Schwartz inequality to \eqref{eq:Differ_f_Pf} yields
	\begin{align*}
		|f(x_j)-\P f(x_j)|
		\leq&~ \Big(\sum_{\ell>n}^{\infty}\sum_{k=1}^{N(d,\ell)}(1+\ell)^{2s}|\widehat{f}_{\ell,k}|^2\Big)^{1/2}\Big(\sum_{\ell>n}^{\infty}\sum_{k=1}^{N(d,\ell)}(1+\ell)^{-2s}|Y_{\ell,k}(x_j)|^2\Big)^{1/2}\\
		\leq &~\|f\|_{H^{
				s}(\S^d)}\Big(\sum_{\ell>n}^{\infty}(1+\ell)^{-2s}\frac{N(d,\ell)}{\omega_d}\Big)^{1/2} \\
		\leq&~ C\omega_d^{-1/2}\|f\|_{H^{
				s}(\S^d)}\Big(\sum_{\ell>n}^{\infty}\ell^{-2s+d-1}\Big)^{1/2}\\
		\leq&~ C\omega_d^{-1/2}n^{-s+d/2}\|f\|_{H^{
				s}(\S^d)},
	\end{align*}
	where we have used that $N(d,\ell)=\mathcal{O}(\ell^{d-1})$ for large $\ell$.
\end{proof}

We also need the following Marcinkiewicz-Zygmund
inequalities for spherical polynomials \cite{dai_2006JFA_multivariate,dai_2006PAMS_generalized}. 
\begin{lemma}[\cite{dai_2006JFA_multivariate}]\label{lem:MZ_Inequality}
	Let $\tau\in\mathbb{N}$ and $X$ be a quasi-uniform set. Suppose that $\{(w_j,x_j)\}$ is a positive quadrature rule on $\S^d$ that is exact for polynomials of degree less than $\tau$. For any $p(x)\in\Pi_{\tau}(\S^d) $ and $0<q<\infty$, there holds
	\begin{equation}
		\begin{aligned}
			c_1\|p\|^q_{L_q(\S^d)}\leq \sum_{j}w_j|p(x_j)|^q\leq c_2\|p\|^q_{L_q(\S^d)}.
		\end{aligned}
	\end{equation}
\end{lemma}

\begin{proof}[\textbf{Proof of Theorem \ref{thm:IntegralError}}]
	
	We split the discretization error into four parts:
	\begin{equation}\label{eq:integralFourParts}
		\begin{aligned}
			&~\big|f*\scaleKer-\sum w_j f(x_j)\scaleKer(x_j\cdot x)\big|\\
			\leq &~\big|f*\scaleKer- \P f*\scaleKer\big|+\big|\P f*\scaleKer-\P f*\P\scaleKer\big|\\
			&~~+|\P f*\P\scaleKer-\sum w_j \P f(x_j)\P\scaleKer(x_j\cdot x)|\\
			&~~+\Big|\sum w_j \P f(x_j)\P\scaleKer(x_j\cdot x)-\sum w_j f(x_j)\scaleKer(x_j\cdot x)\Big|.
		\end{aligned}  
	\end{equation}
	Applying the convolution-type of Young's inequality to the first term and using Lemma \ref{lem:L2Proj} leads to
	\begin{equation}\label{eq:err1}
		\begin{aligned}
			\|f*\scaleKer- \P f*\scaleKer\|_{L_2(\S^d)}
			\leq \|f-\P f\|_{L_2(\S^d)}\|\scaleKer\|_{L_1(\S^d)}
			\leq Cn^{-s}\|f\|_{H^{s}(\S^d)}.
		\end{aligned}
	\end{equation}
	The second and third terms on the right-hand side of \eqref{eq:integralFourParts} vanish because (i) $\P \scaleKer$ is the $L_2$ projection onto $\Pi_n(\S^d)$ and $\P f \in \Pi_n(\S^d)$; (ii) The quadrature rule $\QuadSet$ is of order $2n.$
	Finally, to bound the fourth term in \eqref{eq:integralFourParts}, we rewrite it as
	\begin{equation}\label{eq:err2}
		\begin{aligned}
			&~\Big|\sum w_j \P f(x_j)\P\scaleKer(x_j\cdot x)-\sum w_j f(x_j)\scaleKer(x_j\cdot x)\Big|\\
			=&~\Big|\sum w_j \big(\P f(x_j)-f(x_j)\big)\P\scaleKer(x_j\cdot x)\Big|+\Big|\sum w_j f(x_j)\big(\P\scaleKer(x_j\cdot x)-\scaleKer(x_j\cdot x)\big)\Big|\\
			:=&~\mathcal{E}_{X,f}+\mathcal{E}_{X,\scaleKer}.
		\end{aligned}
	\end{equation}
	The first sum above can be estimated by using the Marcinkiewicz-Zygmund inequality in Lemma \ref{lem:MZ_Inequality} for $q=1$,
	\begin{equation}
		\begin{aligned}
			\mathcal{E}_{X,f}=&~\left|\sum w_j \big(\P f(x_j)-f(x_j)\big)\P\scaleKer(x_j\cdot x)\right|\\
			\leq&~\|\P f-f\|_{L_{\infty}(\S^d)}~\sum\big| w_j\P\scaleKer(x_j\cdot x)\big|\\
			\leq &~ C\|\P f-f\|_{L_{\infty}(\S^d)} \|\P\scaleKer\|_{L_1(\S^d)}.
		\end{aligned}
	\end{equation}
	Based on the $L_2$ error estimate in Lemma \ref{lem:L2Proj}, we can obtain
	\begin{equation}\label{eq:estimateL1}
		\begin{aligned}
			\|\P\scaleKer-\scaleKer\|_{L_1(\S^d)}
			\leq C\|\P\scaleKer-\scaleKer\|_{L_2(\S^d)}
			\leq Cn^{-s}\rho^{-s}
		\end{aligned}
	\end{equation}
	and $\|\P\scaleKer\|_{L_1(\S^d)}\leq \|\P\scaleKer-\scaleKer\|_{L_1(\S^d)}+\|\scaleKer\|_{L_1(\S^d)}$. This together with Lemma \ref{lem:ProjPointWiseErr} leads to
	\begin{equation}\label{eq:err3}
		\begin{aligned}
			\mathcal{E}_{X,f}
			\leq Cn^{-s+d/2}\|f\|_{H^{s}(\S^d)}\Big(n^{-s}\rho^{-s}+\|\scaleKer\|_{L_1(\S^d)}\Big).
		\end{aligned}
	\end{equation}
    Since $n>N_{\rho}=\lfloor \rho^{-1}-1\rfloor$, it follows that $n^{-s}\rho^{-s}\leq 2^s$. From the definition of $\varphi_\rho$ in \eqref{eq:DefiSphKernel}, we have $\|\varphi_\rho\|_{L_1(\S^d)}=1$. Thus, \eqref{eq:err3} simplifies to
    \begin{equation*}
        \mathcal{E}_{X,f}\leq Cn^{-s+d/2}\|f\|_{H^{s}(\S^d)}.
    \end{equation*}
	For the second term, we use Cauchy-Schwartz inequality to obtain
	\begin{equation}\label{eq:err4}
		\begin{aligned}
	\|\mathcal{E}_{X,\scaleKer}\|_{L_2(\S^d)}^2=&~\int_{\S^d}\Big|\sum_j w_j f(x_j)\big(\P\scaleKer(x_j\cdot x)-\scaleKer(x_j\cdot x)\big)\Big|^2\d\mu(x)\\
    \leq &~\int_{\S^d}\Big(\sum_jw_j|f(x_j)|^2\Big) \Big(\sum_jw_j\big|\P\scaleKer(x_j\cdot x)-\scaleKer(x_j\cdot x)\big|^2\Big)\d\mu(x)\\
     \leq &~\|f\|_{L_{\infty}(\S^d)}^2\big(\sum_jw_j\big)\cdot\sum_jw_j\int_{\S^d} \big|\P\scaleKer(x_j\cdot x)-\scaleKer(x_j\cdot x)\big|^2\d\mu(x)\\
			\leq &~\|f\|_{L_{\infty}(\S^d)}^2\|\P\scaleKer-\scaleKer\|_{L_2(\S^d)}^2\big(\sum_jw_j\big)^2\\
			\leq &~Cn^{-2s}\rho^{-2s}\|f\|_{L_{\infty}(\S^d)}^2,
		\end{aligned}
	\end{equation}
    where in the second-to-last inequality, we utilized the rotational invariance of the kernel $\scaleKer$, while in the final inequality, we applied the boundedness of $\sum_jw_j$ along with the estimate from \eqref{eq:estimateL1}.

	Combining inequalities \eqref{eq:err1}, \eqref{eq:err3} and \eqref{eq:err4}, we have 
	\begin{equation*}
		\|f*\varphi_{\rho}-\ConvPhi^Xf\|_{L_2}\leq C\Big(n^{-s}\|f\|_{H^{s}(\S^d)}+n^{-s+d/2}\|f\|_{H^{s}(\S^d)}+n^{-s}\rho^{-s}\|f\|_{L_{\infty}(\S^d)}\Big),
	\end{equation*}
    which completes the proof of Theorem \ref{thm:IntegralError} by using $\|f\|_{L_{\infty}(\S^d)}\leq C\|f\|_{H^{s}(\S^d)}$ for $s>d/2$.
\end{proof}

We summarize the results of Theorem \ref{thm:SphConvConvergence} and Theorem  \ref{thm:IntegralError} in the following theorem and corollary.

\begin{theorem}\label{thm:SphQuasiInterpErr}
	Let a  quasi-interpolant be given as in \eqref{eq:Spherical_quasi_interp} with a scaled zonal kernel $\scaleKer$ that satisfies Assumption \ref{Assump1:Kernel} and Assumption \ref{Assump2:Kernel} with $ s>d/2$. Suppose that $X$ is a quasi-uniform set and $\mathcal{Q}_X$ is a positive quadrature rule of order $2n$ with $n > N_\rho$.  Then there exists a constant $C$ independent of $n$ and $\rho$ such that,
	\begin{align}\label{eq:estimateQuasiInterp}
		\|f-\ConvPhi^Xf\|_{L_2(\S^d)}\leq &~C\big(\rho^{s}+n^{-s}+n^{-s+d/2}+n^{-s}\rho^{-s}\big)\|f\|_{H^{s}(\S^d)}.
	\end{align}
\end{theorem}


\begin{corollary}\label{Corol:1}
	Under the assumptions of Theorem \ref{thm:SphQuasiInterpErr}, if we choose the scaled parameter $\rho=\mathcal{O}(n^{-1/2})$, then the quasi-interpolation operator \eqref{eq:Spherical_quasi_interp} has the following error estimate
	$$\|f-\ConvPhi^X f\|_{L_2(\S^d)}\leq Cn^{-\min\{s/2,~ s-d/2\}}\|f\|_{H^s(\S^d)}.$$
\end{corollary}
\begin{proof}
   Since $\rho<1$, we have $n^{-s}\rho^{-s}>n^{-s}$. Thus, letting $\rho^{s}= \mathcal{O}(n^{-s}\rho^{-s})$ in \eqref{eq:estimateQuasiInterp} completes the proof.
\end{proof}

\section{Numerical examples}
\label{sec:NumerExamp}
In this section, we present several numeral examples to demonstrate theoretical assertions  for the proposed spherical quasi-interpolation method. Our examples employ the maximum determinant (MD) nodes on the sphere $\mathbb{S}^2$, as introduced by Womersley and Sloan \cite{womersley_sloan_spherepoints}. The MD nodes are known for their excellent conditioning, and the associated cubature weights are numerically positive. The weights computed in \cite{womersley_sloan_spherepoints} ensure that the cubature rule is exact for all polynomials of degree less than $n$, where $|X| = (n+1)^2$ is the corresponding number of MD nodes used. The numerical $L_2$ error is estimated on a set of $|Y|=36864$ evaluation nodes. We use a discrete $L_2(\mathbb{S}^d)$-norm based on a quadrature rule defined at $\{y_j\}_{j=1}^{|Y|}$, a common practice in the literature \cite{fuselier2013high,kunemund2019high}:
$\|f\|_{L_2(\mathbb{S}^d)}^2 \approx \sum_{j=1}^{|Y|} w_j[f(y_j)]^2$.

\subsection{Convergence test} To evaluate the performance of the proposed quasi-interpolation method, we consider three distinct classes of functions on the sphere: 

\textbf{Spherical harmonic function \cite{fuselier2013high}}:
\begin{equation}
	f_1(x)=Y_{6,4}(x), ~x\in\S^d.
\end{equation}

\textbf{Gaussian function with random centers \cite{calhoun2010finite,fuselier2013high}}:
\begin{equation}\label{eq:randGaussian}
	f_2(x) = \sum_{k=1}^{\Lambda}\exp(-10\arccos(\xi_k\cdot x)),
\end{equation}
where $\xi_k$, $k=1,2,\ldots,\Lambda$ are randomly placed points on the surface of the sphere.

\textbf{Finitely smooth function \cite{franz-wendland2023multilevel}:}  
\begin{equation}\label{eq:finiteSmoothFunc}
	f_3(x)= \Big(1-\sqrt{2-2\cos(x)}\Big)^2_{+}.
\end{equation}

In the first example, we employ the Poisson kernel to construct a quasi-interpolation scheme for approximating the three classes of functions on the sphere. Table \ref{tab:convOrder_poisson} shows the numerical errors and convergence orders for $n=10, 20, 40, 80, 160$. The results demonstrate that the Poisson kernel quasi-interpolation achieves a convergence order of approximately one-half, which aligns with the theoretical rate of $\mathcal{O}(n^{-1/2})$ in Corollary \ref{Corol:1}, as it only satisfies Assumption \ref{Assump1:Kernel} with $s=1$.
Additionally, we test the performance of Gaussian kernels and compactly-supported kernels of varying orders in approximating the spherical harmonic function $Y_{6,4}$.  The numerical errors and convergence rates are presented in Table \ref{tab:convOrder_sphHarm} for restricted Gaussian kernels and in Table \ref{tab:convOrder_sphHarm_cswu} for compactly-supported positive definite (CSPD) functions \cite[Example 3]{gao-zhou2020multiscale} satisfying different orders of Assumption \ref{Assump1:Kernel}. It is evident that when choosing $\rho=\mathcal{O}(n^{-1/2})$, the approximation orders are consistent with our theoretical error bounds.

\begin{table}[htbp]
	\centering
	\caption{Convergence results of spherical quasi-interpolation method for approximating $f_1(x)$, $f_2(x)$ and $f_3(x)$ by using \textbf{Poisson kernel} with shape parameter $\rho=0.4n^{-1/2}$.}
	\label{tab:convOrder_poisson}
	\begin{tabular}{*{7}{c}}
		\toprule
		$n$& $f_1(x)$ &rate & $f_2(x)$ &rate & $f_3(x)$ & rate  \\
		\midrule
		$10$&  2.3591e-01   & -& 2.5379e-01 & -&1.2884e-01 & - \\
		$20$&  1.7432e-01   & 0.44& 1.6694e-01 & 0.60&7.7334e-02 & 0.74 \\
		
		$40$&   1.2963e-02 &0.43 & 1.1684e-01 & 0.51 & 5.4609e-02 & 0.50 \\
		
		$80$&  9.5783e-02  & 0.44& 8.5147e-02 & 0.46& 3.9611e-02 &  0.46\\
		$160$&  6.9957e-02 & 0.45 & 6.1975e-02& 0.46& 2.8804e-02  & 0.46 \\
		\bottomrule
	\end{tabular}
\end{table}

\begin{table}[htbp]
	\centering
	\caption{Convergence results of spherical quasi-interpolation method for approximating spherical harmonics $Y_{6,4}$ by using various orders of \textbf{restricted Gaussian kernel} ($s=2,4,6$) with shape parameter $\rho=0.4n^{-1/2}, 0.7n^{-1/2}, 1n^{-1/2}$.}
	\label{tab:convOrder_sphHarm}
	\begin{tabular}{*{7}{c}}
		\toprule
		$n$& $s=2$ &rate & $s=4$ &rate & $s=6$ & rate  \\
		\midrule
		$10$&  3.5633e-01   & -& 2.6752e-01 & -&1.7896e-01 & - \\
		$20$&  1.5921e-01   & 1.16& 4.9910e-02 & 2.42&2.9474e-02 & 2.60 \\
		
		$40$&   8.0797e-02 &0.98 & 9.9009e-03 & 2.34 & 3.6412e-03 & 3.02 \\
		
		$80$&  4.1170e-02  & 0.97& 2.5158e-03 & 1.98& 4.9753e-04 &  2.87\\
		$160$&  2.0791e-02 & 0.99 & 6.4259e-04  & 1.97& 6.5216e-05  & 2.93 \\
		\bottomrule
	\end{tabular}
\end{table}

\begin{table}[htbp]
	\centering
	\caption{Convergence results of spherical quasi-interpolation method for approximating spherical harmonics $Y_{6,4}$ by using various orders of \textbf{restricted compactly-supported kernel} ($s=2,4,6$) with shape parameter $\rho=1.5n^{-1/2}, 3n^{-1/2}, 4n^{-1/2}$.}
	\label{tab:convOrder_sphHarm_cswu}
	\begin{tabular}{*{7}{c}}
		\toprule
		$n$& $s=2$ &rate & $s=4$ &rate & $s=6$ & rate  \\
		\midrule
		$10$&  3.7864e-01   & -& 2.3209e-01 & -&1.6278e-01 & - \\
		$20$&  1.2889e-01   & 1.55& 4.3025e-02 & 2.43&2.3111e-02 & 2.82 \\
		
		$40$&   6.1108e-02 &1.08 & 8.4685e-03 & 2.35 &  1.7921e-03 & 3.69 \\
		
		$80$&  3.0847e-02  & 0.99& 2.1855e-03 & 1.95& 2.3795e-04 &  2.91\\
		$160$&  1.5529e-02 & 0.99 & 5.5594e-04  & 1.98& 3.0671e-05  & 2.96 \\
		\bottomrule
	\end{tabular}
\end{table}

\begin{figure}[htbp]
	\centering
	\begin{tikzpicture}[scale=0.8]
		\begin{loglogaxis}[
			grid=both,
			grid style = dotted,
			mark size = 1.5pt,
			xmin = 1, xmax = 10^3,
			ymin = 0.5*10^(-5),ymax = 10^(-0),
			xlabel={$n$},
			title={$f_2(x)$, Gaussian kernel},
			legend cell align = {left},
			legend pos = south west,
			legend entries={$s=2$,$s=4$,$s=6$}]
			\addplot  +[ultra thick,dashed, mark = square*,mark options={solid},red] table [x=n,y=L2err_m2] {Data_sphereQI/L2err_sphere_gauss_randGaussian.dat};
			\addplot  +[ultra thick,  dotted, mark = *, color = green] table [x=n,y=L2err_m4] {Data_sphereQI/L2err_sphere_gauss_randGaussian.dat};
			\addplot  +[ultra thick,mark = triangle*, blue] table [x=n,y=L2err_m6] {Data_sphereQI/L2err_sphere_gauss_randGaussian.dat};
			\addplot[dashed, color = black] coordinates {(40,0.7*1e-1) (160,0.7*1e-1/4) };
			\addplot[dashed, color = black] coordinates {(40,0.9*1e-2) (160,0.9*1e-2/16)};
			\addplot[dashed, color = black] coordinates {(40, 0.8*1e-3) (160,0.8*1e-3/64)};
		\end{loglogaxis}
		\node[rectangle] at (5.4,4) {\scriptsize $n^{-1}$};
		\node[rectangle] at (5.4,2.5) {\scriptsize $n^{-2}$};
		\node[rectangle] at (5.4,0.5) {\scriptsize $n^{-3}$};
	\end{tikzpicture}
	\hspace{0.5cm}
	\begin{tikzpicture}[scale=0.8]
		\begin{loglogaxis}[
			grid=both,
			grid style = dotted,
			mark size = 1.5pt,
			xmin = 1, xmax = 10^3,
			ymin = 0.5*10^(-5),ymax = 10^(-0),
			xlabel={$n$},
			title={$f_2(x)$, Compactly-supported kernel},
			legend cell align = {left},
			legend pos = south west,
			legend entries={$s=2$,$s=4$,$s=6$}]
			\addplot  +[ultra thick,dashed, mark = square*,mark options={solid},red] table [x=n,y=L2err_m2] {Data_sphereQI/L2err_sphere_cswu_randGauss.dat};
			\addplot  +[ultra thick,  dotted, mark = *, color = green] table [x=n,y=L2err_m4] {Data_sphereQI/L2err_sphere_cswu_randGauss.dat};
			\addplot  +[ultra thick,mark = triangle*, blue] table [x=n,y=L2err_m6] {Data_sphereQI/L2err_sphere_cswu_randGauss.dat};
			\addplot[dashed, color = black] coordinates {(40,0.7*1e-1) (160,0.7*1e-1/4) };
			\addplot[dashed, color = black] coordinates {(40,0.9*1e-2) (160,0.9*1e-2/16)};
			\addplot[dashed, color = black] coordinates {(40, 0.8*1e-3) (160,0.8*1e-3/64)};
		\end{loglogaxis}
		\node[rectangle] at (5.4,4) {\scriptsize $n^{-1}$};
		\node[rectangle] at (5.4,2.5) {\scriptsize $n^{-2}$};
		\node[rectangle] at (5.4,0.5) {\scriptsize $n^{-3}$};
	\end{tikzpicture}
	\vspace{-1pt}
	\caption{Approximation errors and convergence orders for the target function \eqref{eq:randGaussian} using spherical quasi-interpolation with various numbers of MD nodes and kernels satisfying Assumption \ref{Assump1:Kernel} with orders $s=2,4,6$. The left figure corresponds to restricted Gaussian kernels, and the right corresponds to compactly-supported kernels.}
	\label{fig.L2conv_randGaussian}
\end{figure}
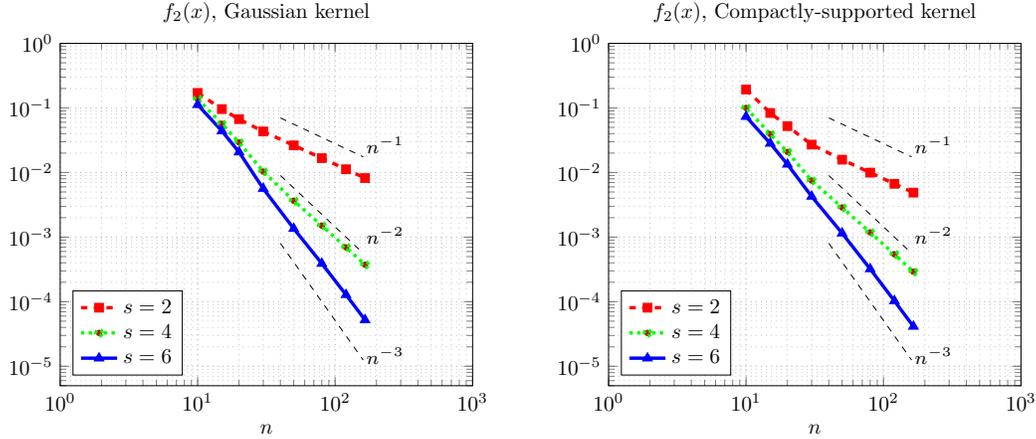

\begin{figure}[htbp]
	\centering
	\begin{tabular}{ccccc}
		\begin{overpic}[width=0.45\textwidth,trim= 0 0 0 0, clip=true,tics=10]{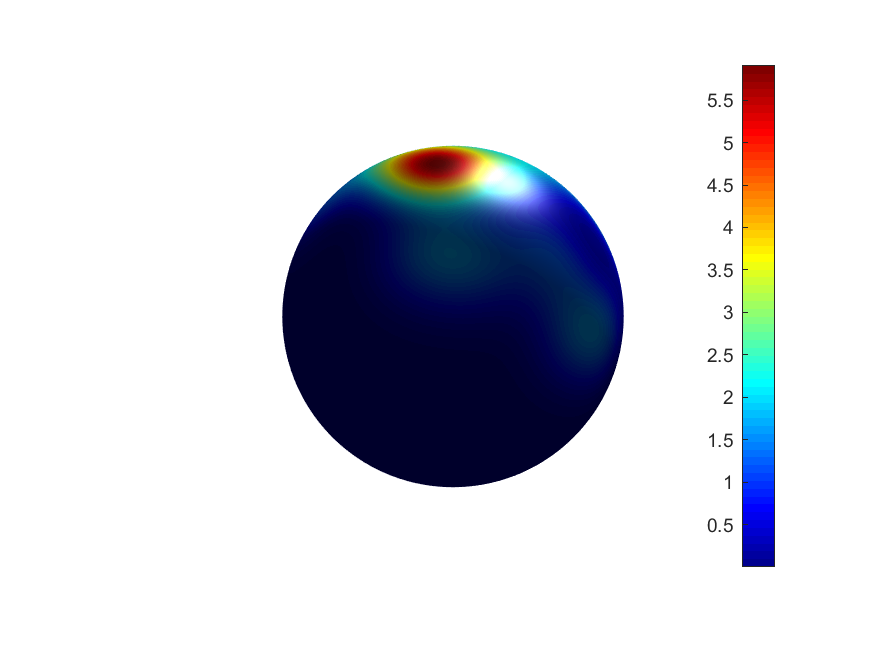}
		\put (24,-2) {(a) Numerical solution}
	\end{overpic}
	&
	\begin{overpic}[width=0.45\textwidth,trim= 0 0 0 0, clip=true,tics=10]{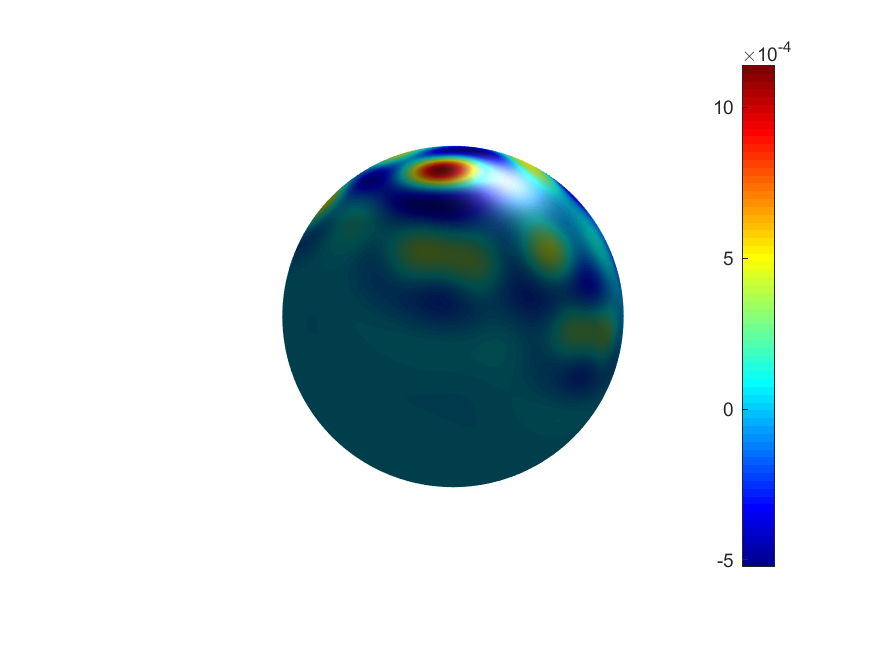}
	\put (24,-2) {(b) Approximation error}
\end{overpic}
\end{tabular}
\caption{Numerical solution and approximation error of the spherical quasi-interpolation for approximating \eqref{eq:randGaussian} using a restricted Gaussian kernel satisfying Assumption \ref{Assump1:Kernel} with $s=6$ and $|X|=101^2$ MD nodes.} 
\label{fig.NumerSolu}
\end{figure}

The second example considers a test problem presented in \cite{calhoun2010finite,fuselier2013high} to demonstrate the approximation property on the sphere. The target function is given by equation \eqref{eq:randGaussian}, where each term is a Gaussian centered at $\xi_k$, with the distance measured using the geodesic distance. The solution is clearly $C^{\infty}(\mathbb{S}^2)$. In our simulation, we choose $\Lambda=23$ as in \cite{fuselier2013high}. The numerical errors and convergence results are presented in Figure \ref{fig.L2conv_randGaussian}, which again verifies the theoretical results. The numerical solution and pointwise approximation error are illustrated in Figure \ref{fig.NumerSolu}.

In our last example, the test function $f_3(x)$ given by equation \eqref{eq:finiteSmoothFunc} is not particularly smooth. This implies that it does not satisfy the required conditions for our error analysis for larger values of $s$. The results are presented in Figure \ref{fig.L2conv_finiteSmooth}, which demonstrates first-order of convergence for both kernels with all values of $s$.

\definecolor{ao(English)}{rgb}{0.0,0.5,0.0}

\begin{figure}[htbp]
\centering
\begin{tikzpicture}[scale=0.8]
\begin{loglogaxis}[
	grid=both,
	grid style = dotted,
	mark size = 1.5pt,
	xmin = 1, xmax = 10^3,
	ymin = 0.5*10^(-3),ymax = 10^(-0),
	xlabel={$n$},
	title={$f_3(x)$, Gaussian kernel},
	legend cell align = {left},
	legend pos = south west,
	legend entries={$s=2$,$s=4$,$s=6$}]
	\addplot  +[ultra thick,dashed, mark = square*,mark options={solid},red] table [x=n,y=L2err_m2] {Data_sphereQI/L2err_sphere_gauss_finitesmooth.dat};
	\addplot  +[ultra thick,  dotted, mark = *, color = green] table [x=n,y=L2err_m4] {Data_sphereQI/L2err_sphere_gauss_finitesmooth.dat};
	\addplot  +[ultra thick,mark = triangle*, blue] table [x=n,y=L2err_m6] {Data_sphereQI/L2err_sphere_gauss_finitesmooth.dat};
	\addplot[dashed, color = black] coordinates {(40,0.7*1e-1) (160,0.7*1e-1/4) };
\end{loglogaxis}
\node[rectangle] at (5.4,3) {\scriptsize $n^{-1}$};
\end{tikzpicture}
\hspace{0.5cm}
\begin{tikzpicture}[scale=0.8]
\begin{loglogaxis}[
	grid=both,
	grid style = dotted,
	mark size = 1.5pt,
	xmin = 1, xmax = 10^3,
	ymin = 0.5*10^(-3),ymax = 10^(-0),
	xlabel={$n$},
	title={$f_3(x)$, Compactly-supported function},
	legend cell align = {left},
	legend pos = south west,
	legend entries={$s=2$,$s=4$,$s=6$}]
	\addplot  +[ultra thick,dashed, mark = square*,mark options={solid},red] table [x=n,y=L2err_m2] {Data_sphereQI/L2err_sphere_cswu_finitesmooth.dat};
	\addplot  +[ultra thick,  dotted, mark = *, color = green] table [x=n,y=L2err_m4] {Data_sphereQI/L2err_sphere_cswu_finitesmooth.dat};
	\addplot  +[ultra thick,mark = triangle*, blue] table [x=n,y=L2err_m6] {Data_sphereQI/L2err_sphere_cswu_finitesmooth.dat};
	\addplot[dashed, color = black] coordinates {(40,0.6*1e-1) (160,0.6*1e-1/4) };
\end{loglogaxis}
\node[rectangle] at (5.4,2.8) {\scriptsize $n^{-1}$};
\end{tikzpicture}
\vspace{-1pt}
\caption{Approximation errors and convergence orders for the target function \eqref{eq:finiteSmoothFunc} using spherical quasi-interpolation with various numbers of maximum determinant nodes and kernels satisfying Assumption \ref{Assump1:Kernel} with orders $s=2,4,6$, including restricted Gaussian kernels and CSPD functions.}
\label{fig.L2conv_finiteSmooth}
\end{figure}
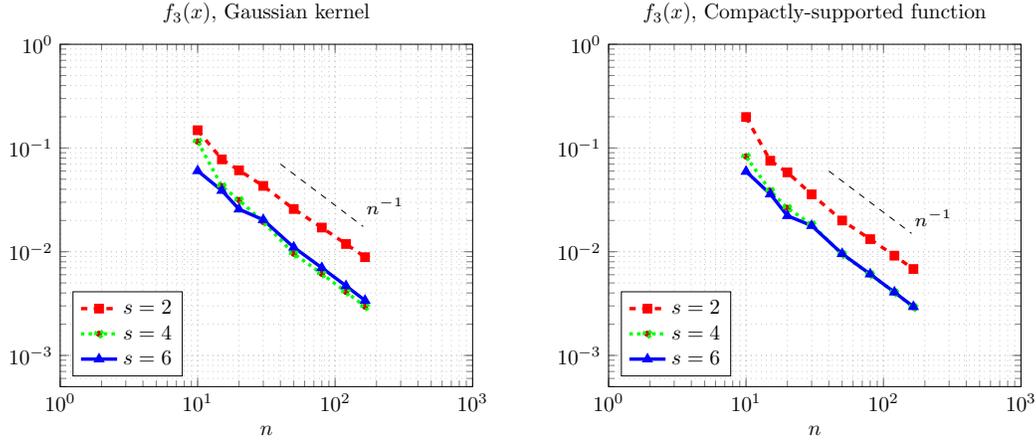

\begin{figure}[htbp]
\centering
\begin{tikzpicture}[scale=0.81]
\begin{loglogaxis}[
	grid=both,
	grid style = dotted,
	mark size = 1.5pt,
	ymin = 10^(-4),ymax = 2*10^0,
	xlabel={$|X|$},
	title={Spherical hyperinterpolation},
	legend cell align = {left},
	legend pos = south west,
	legend columns=3,
	legend style = {font=\footnotesize},
	legend entries={{ $\delta=0.001$},{$\delta=0.01$}, { $\delta=0.1$},{ $\delta=0.3$},{ $\delta=0.5$}}]
	\addplot  +[ultra thick,  color = r1] table [x=N,y=L2err1] {Data_sphereQI/L2err_hyper_noise.dat};
	\addplot  +[ultra thick,   color = Orange] table [x=N,y=L2err2] {Data_sphereQI/L2err_hyper_noise.dat};
	\addplot  +[ultra thick,  color = g1] table [x=N,y=L2err3] {Data_sphereQI/L2err_hyper_noise.dat};
	\addplot  +[ultra thick,  color = b1] table [x=N,y=L2err4] {Data_sphereQI/L2err_hyper_noise.dat};
	\addplot  +[ultra thick,  color = JungleGreen] table [x=N,y=L2err5] {Data_sphereQI/L2err_hyper_noise.dat};
\end{loglogaxis}
\end{tikzpicture}
\hspace{0.7cm}
\begin{tikzpicture}[scale=0.8]
\begin{loglogaxis}[
	grid=both,
	grid style = dotted,
	mark size = 1.5pt,
	ymin = 10^(-4),ymax = 2*10^0,
	xlabel={$|X|$},
	title={Spherical quasi-interpolation},
	legend cell align = {left},
	legend pos = south west,
	legend style = {font=\footnotesize},
	legend columns=3,
	legend entries={{ $\delta=0.001$},{$\delta=0.01$}, { $\delta=0.1$},{ $\delta=0.3$},{ $\delta=0.5$}}]
	\addplot  +[ultra thick,  color = r1] table [x=N,y=L2err1] {Data_sphereQI/L2err_sphQI_noise.dat};
	\addplot  +[ultra thick,   color = Orange] table [x=N,y=L2err2] {Data_sphereQI/L2err_sphQI_noise.dat};
	\addplot  +[ultra thick, color = g1] table [x=N,y=L2err3] {Data_sphereQI/L2err_sphQI_noise.dat};
	\addplot  +[ultra thick,  color = b1] table [x=N,y=L2err4] {Data_sphereQI/L2err_sphQI_noise.dat};
	\addplot  +[ultra thick,  color = JungleGreen] table [x=N,y=L2err5] {Data_sphereQI/L2err_sphQI_noise.dat};
\end{loglogaxis}
\end{tikzpicture}
\vspace{-1pt}
\caption{RMSE error performance of spherical quasi-interpolation versus spherical hyperinterpolation for approximating a noisy target function. The approximation error is measured for different noise levels $\delta=0.001,0.01,0.1,0.3,0.5$ and various MD nodes $|X|=5^2,7^2,\ldots, 161^2$.}
\label{fig.Compar_noisydata}
\end{figure}
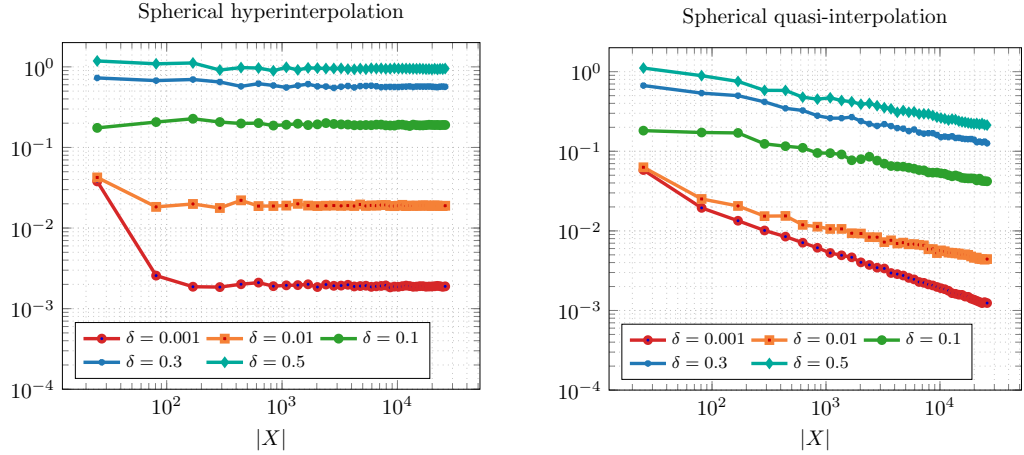

\pgfplotstableread[row sep=\\,col sep=&]{
N & hyperT & qiT  \\
3721     & 1.8  & 1.9    \\
6561     & 3.3 & 3.0    \\
10201    & 5.4 & 4.4  \\
14641   & 8.4 & 5.8  \\
19881   & 11.2  & 6.9  \\
25921      & 17.3  & 9.3 \\
}\mydata

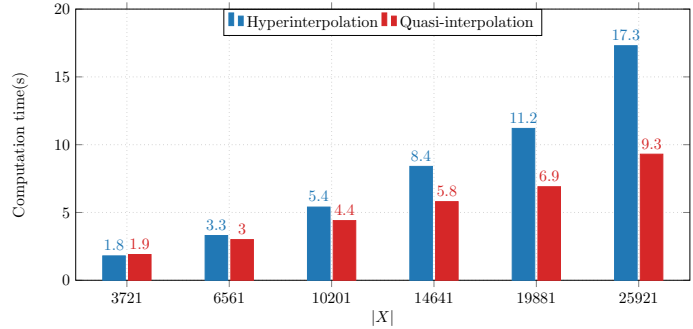
\begin{figure}[htbp]
\centering
\begin{tikzpicture}[scale=0.6]
\begin{axis}[
	ybar,
	bar width=.5cm,
	width=\textwidth,
	height=.5\textwidth,
	legend style={at={(0.5,1)},
		anchor=north,legend columns=-1},
	symbolic x coords={3721,6561,10201,14641,19881,25921},
	xlabel={$|X|$},
	xtick=data,
	grid = major,
	grid style = dotted,
	nodes near coords,
	nodes near coords align={vertical},
	ymin=0,ymax=20,
	ylabel={Computation time(s)},
	]
	\addplot+ [color = b1] table[x=N,y=hyperT]{\mydata};
	\addplot+ [color = r1] table[x=N,y=qiT]{\mydata};
	\legend{Hyperinterpolation, Quasi-interpolation}
\end{axis}
\end{tikzpicture}
\vspace{-1pt}
\caption{Computation time comparison between spherical quasi-interpolation (red) and spherical hyperinterpolation (blue) for different numbers of MD nodes ($|X|=61^2,81^2,101^2,121^2,141^2,161^2$). }
\label{fig.Compar_CpuTime}
\end{figure}

\subsection{Simulation for noisy data}
In this test, we aim to evaluate and compare the effectiveness of two distinct methods for approximating functions with additive noise on a spherical domain: the proposed spherical quasi-interpolation method using a Gaussian kernel and the spherical hyperinterpolation technique developed by Sloan in \cite{sloan_1995JAT_polynomial}. Our primary goal is to assess the accuracy and computational efficiency of these methods under varying conditions. To this end, we conduct simulations using different numbers of MD nodes, ranging from $5^2$ to $161^2$, and analyze their performance in terms of approximation error and computation time. The target function used in our study is defined as:
\begin{equation}
f(x) = \sum_{j=1}^6\phi_3(\|x-z_j\|),
\end{equation}
where $\phi_3$ is a normalized Wendland function as defined in \cite{wendland2004scattered},  and the points $z_1,z_2=[\pm 1,0,0]^T$, $z_3,z_4=[ 0,\pm1,0]^T$ and $z_5,z_6=[ 0,0, \pm1]^T$. This function was previously used to test the behavior of hyperinterpolation in \cite{an-BIT2022-quadrature}. To simulate noisy data, we add Gaussian-distributed noise $\mathcal{N}(0,\delta^2)$ at five different levels: $\delta=0.001,0.01,0.1,0.3,0.5$ to the above function as sampling data. 

The numerical results are presented in Figure \ref{fig.Compar_noisydata}, where the RMSE error is measured on a very fine node set ($N=50000$) to assess the approximation error. Each simulation is repeated $30$ times, and the average results are recorded. Figure \ref{fig.Compar_noisydata} shows that hyperinterpolation converges for small noise levels but deteriorates as the noise increases. In contrast, our spherical quasi-interpolation method shows a consistently decreasing RMSE as the number of sampling nodes grows. Moreover, Figure \ref{fig.Compar_CpuTime} compares the computation time for spherical hyperinterpolation and the proposed spherical quasi-interpolation method with various numbers of MD nodes ($|X|=61^2,81^2,101^2,121^2,141^2,161^2$). The results indicate that as the number of sampling nodes increases, the computation time for the proposed quasi-interpolation method is significantly less than that of the hyperinterpolation method. For $|X|=161^2$, the computation time of the hyperinterpolation method is roughly double that of the proposed method. This demonstrates that our method is both robust and accurate for approximating noisy data.

\subsection*{Acknowledgements}
The work of the first author was supported  by  NSFC (No.12101310), NSF of Jiangsu Province (No.BK20210315), and the Fundamental Research Funds for the Central Universities (No.30923010912), The work of the second author was supported by NSFC (No.12271002).



\bibliographystyle{plain}
\bibliography{periodic_kernel}

\end{document}